\documentclass[number,citesort,seceqn,dvips]{arxbj}

\aid{0}
\volume{17}
\issue{3}
\pubyear{2011}
\firstpage{1095}
\lastpage{1125}
\doi{10.3150/10-BEJ305}

\makeatletter

\newtheorem{theorem}{Theorem}[section]

\newtheorem{corollary}[theorem]{Corollary}
\newtheorem{proposition}[theorem]{Proposition}
\newproclaim{definition}[theorem]{Definition}
\newremark{remark}[theorem]{Remark}

\newremark{example}[theorem]{Example}

\makeatother

\begin{document}
\begin{frontmatter}

\title{Multivariate Jacobi and Laguerre polynomials,
infinite-dimensional extensions, and
their probabilistic connections with multivariate Hahn and Meixner polynomials}
\runtitle{Multivariate Jacobi}

\begin{aug}
\author[a]{\fnms{Robert C.} \snm{Griffiths}\thanksref{a}\ead[label=e1]{griff@stats.ox.ac.uk}}
\and
\author[b]{\fnms{Dario} \snm{Span\`{o}}\corref{}\thanksref{b}\ead[label=e2]{d.spano@warwick.ac.uk}}
\runauthor{R.C. Griffiths and D. Span\`{o}}
\address[a]{Department of Statistics, 1 South Parks Road Oxford OX1
3TG, UK.\\ \printead{e1}}
\address[b]{Department of Statistics, University of Warwick, Coventry
CV4 7AL, UK.\\ \printead{e2}}
\end{aug}

\received{\smonth{11} \syear{2009}}
\revised{\smonth{7} \syear{2010}}

%
\begin{abstract}
Multivariate versions of classical orthogonal polynomials such as
Jacobi, Hahn, Laguerre and Meixner are reviewed and their connection
explored by adopting a probabilistic approach. Hahn and Meixner
polynomials are interpreted as posterior mixtures of Jacobi and
Laguerre polynomials, respectively. By using known properties of
gamma point processes and related transformations, a new
infinite-dimensional version of Jacobi polynomials is constructed
with respect to the size-biased version of the Poisson--Dirichlet
weight measure and to the law of the gamma point process from
which it is derived.
\end{abstract}

%

\end{frontmatter}

\section{Introduction}\label{sec1}

 In this paper we will review multivariate orthogonal
polynomials, complete with respect to weight measures given by the
Dirichlet and Dirichlet-multinomial probability distributions (denoted
respectively as
$D_{\alpha}$ or $\mathit{DM}_\alpha$, $\alpha\in\mathbb{R}_+^d$), that is,
polynomials
$\{G_n\dvtx n\in\mathbb{N}^{d}\}$ satisfying
%
\begin{equation} \label{opdef}
\int G_n G_m
\,\mathrm{d}\mu=\frac{1}{c_m}\delta_{nm},
\qquad n,m\in\mathbb{N}^{d}.
\end{equation}
The polynomials $\{G_n\}$ are known as multivariate Jacobi polynomials if
(\ref{opdef}) is satisfied with $\mu=D_{\alpha}$,
and multivariate Hahn polynomials if $\mu=\mathit{DM}_{\alpha}$.
Here $c_m$
are positive constants. Completeness means that, for every
function $f$ with finite variance (under $\mu$), there is an
expansion
%
\begin{equation}\label{fou}
 f(x)=\sum_{n\in\mathbb{N}^{d}}c_na_nG_n(x),
\end{equation}
where
\[
a_n=\mathbb{E}[f(X)G_n(X)].
\]
Systems of multivariate orthogonal polynomials are not unique, and a
large number of characterizations of $d$-dimensional Jacobi and Hahn
polynomials exist in literature. We will focus on a construction of
Jacobi polynomials, based on a method originally proposed by Koornwinder
\cite{Koo75} that has a strong probabilistic interpretation. Based on
this, we will re-interpret the role of Jacobi polynomials in the
construction of multivariate Hahn and several other well-known classes
of multivariate orthogonal polynomials. In particular, we will
(1) describe multivariate Hahn polynomials
as \emph{posterior} mixtures of Jacobi polynomials, in a sense which
will become precise in Section \ref{sec:hahn}; (2) construct, in Section
\ref{sec:lag}, a new system of multiple Laguerre polynomials,
orthogonal with respect to the product of several gamma probability
distributions with identical scale parameters;
(3)~derive, in Section \ref{sec:meix}, a new class of multiple Meixner
polynomials as posterior mixtures of the
Laguerre polynomials mentioned in (2); (4) obtain polynomials in the
multivariate hypergeometric distribution by taking the parameters in
the Hahn
polynomials to be negative;
(5)~obtain (Section \ref{sec:rank}) asymptotic results as the
dimension $d\rightarrow\infty$ with $|\alpha|:=\sum_{i=1}^{d}\alpha
_i\rightarrow|\theta|>0,$ by considering size-biased Dirichlet
measures.

Furthermore, we will see that an extensive application of Koornwinder's
method leads directly to finding new systems of polynomials, orthogonal
with respect to a wider family of distributions on the infinite
simplex, known in Bayesian nonparametric statistics as the (discrete)
beta-Stacy family \cite{WM97}, a popular member of which is the GEM
distribution (so named after Griffiths, Engen and
McCloskey who introduced it independently) and its two-parameter distribution.

The intricate relationship existing among all the mentioned systems of
polynomials is traditionally described in terms of their
analytic/algebraic expression as (multivariate) basic hypergeometric
series (see, e.g., \cite{Ex76,DX02}). The main advantage of a
probabilistic approach is that it re-expresses most relationships in
terms of random variables, which may be more transparent to
statisticians and probabilists. With this in mind we will begin the
paper with an
introductory summary (Section \ref{sec:distn}) of known facts from the
theory of probability distributions. Section \ref{*} is devoted to
multivariate Jacobi
polynomials, whose structure will be the building block for the
subsequent sections: Multiple Laguerre in Section \ref{sec:lag}, Hahn
in Section \ref{sec:hahn} and Meixner in Section~\ref{sec:meix}.

It is worth observing that the posterior mixture representation of
multivariate Hahn polynomials shown in Proposition \ref{prp:mh} is
obtained without imposing \emph{a priori} any Bernstein--B\'{e}zier
form to the Jacobi polynomials, and nevertheless it agrees with
recent interpretations of Hahn polynomials as Bernstein
coefficients of Jacobi polynomials in such a form \cite{W06,S07},
a result for which a new, more probabilistic proof is
offered in Section \ref{sec:bb1}. In particular, our approach will
make more intuitive the link between the Bernstein--B\'{e}zier
interpretation and the
original formulation proposed decades ago by
Karlin and McGregor
\cite{KMG75}. In terms of applications, understanding such a link
will complete Karlin and McGregor's analysis of some well-known
$d$-type models in population genetics (Section \ref{sec:kmg}). Our
extensions of Sections \ref{sec:rank} and \ref{sec:sblag} open for
possible new infinite-dimensional versions of Karlin and McGregor's
work.

Along the same lines one can
view the Meixner polynomials obtained in Proposition~%
\ref{prp:mmeix} as re-scaled Bernstein
coefficients of our multiple Laguerre polynomials, as shown in Section
\ref{sec:bb2}.

The original motivation for this study was to obtain
some background material that can be used to characterize
bivariate distributions, or transition functions, with fixed
Dirichlet or Dirichlet-multinomial marginals, for which the
following \emph{canonical expansions} are possible:
\begin{eqnarray*}
p(\mathrm{d}x,\mathrm{d}y)=\Biggl\{1+\sum_{n\in\mathbb{Z}^d_+}^{
\infty
}c_n\rho_n
G_n(x)G_n(y)\Biggr\}D_{\alpha}(\mathrm{d}x)D_{\alpha}(\mathrm{d}y),\qquad x,y\in
\Delta_{(d-1)},
\end{eqnarray*}
for appropriate, positive-definite sequences
$\rho_{m}\dvtx m\in\mathbb{N}^{d},$ called the \emph{canonical
correlation coefficients} of the model. Some results on such a
problem are in \cite{GS09} and \cite{GS10}. Other
possible applications in statistics are related to least-squares
approximations and regression. An MCMC (Markov chain Monte Carlo)-Gibbs
sampler use of
orthogonal polynomials is explored, for example, in \cite{SCDK09};
related applications are in \cite{KZ09}.
In this paper, however, we will focus merely on the construction of the
mentioned systems of polynomials.

\section{Distributions on the discrete and continuous simplex}\label{sec:distn}

Throughout the paper we will denote by $|x|$ the total sum of all
components of $x=(x_1,\ldots,x_d)\in\mathbb{R}^d.$ We will also adopt
the notation:
\[
x^\alpha=x_1^{\alpha_1}\cdots x_d^{\alpha_d},\qquad \Gamma
(\alpha)=\prod_{i=1}^{d}\Gamma(\alpha_i)
\]
and
\[
\pmatrix{{|n|}\cr  n}=\frac{|n|!}{\prod_{i=1}^dn_i!}.
\]
For example, the Dirichlet distribution $D_\alpha\dvtx \alpha\in\mathbb
{R}_+^d$ will be written as
\[
D_\alpha(\mathrm{d}x)=\frac{\Gamma(|\alpha|)x^{\alpha-\underline
{1}}}{\Gamma(\alpha)}\mathbb{I}\bigl(x\in\Delta_{(d-1)}\bigr)\,\mathrm{d}x,
\]
where $\underline{1}=(1,1,\ldots,1)$ and, for $d=2,3,\ldots,$
$\Delta_{(d-1)}=\{x\in\mathbb{R}_+^d\dvtx |x|=1\}$.

\subsection{Conditional independence in the Dirichlet distribution}
\subsubsection{Gamma sums}\label{sec:gammasums}

For every $\alpha=(\alpha_1,\ldots,\alpha_d)\in\mathbb{R}^{d}_{+}$
and $\beta>0$, let
$Y=(Y_1,\ldots,Y_{d})$ be a collection of~$d$-inde\-pendent gamma
random variables with parameter, respectively,
$(\alpha_i,\beta).$ The distribution of $Y$ is given by the product
measure
\[
\gamma^{d}_{\alpha,\beta}(\mathrm{d}y)=\frac{y^{\alpha-\underline{1}}
\mathrm{e}^{-{|y|/\beta}}}{\Gamma(\alpha)\beta^{|\alpha|}}
\mathbb{I}(y\in\mathbb{R}_+^d)\, \mathrm{d}y.
\]
Consider the mapping
\[
(Y_1,\ldots,Y_d)\longmapsto (|Y|,X_1,\ldots,X_{d-1}),
\]
where
\[
X_{j}:=\frac{Y_{j}}{|Y|},\qquad j=1,\ldots,d-1,
\]
and set $X_d=1-\sum_{i=1}^{d-1}X_i.$ It is easy to rewrite
\[
\gamma^{d}_{\alpha,\beta}(\mathrm{d}y)=\gamma^1_{|\alpha|,\beta
}(\mathrm{d}|y|)D_{\alpha}(\mathrm{d}x),
\]
that is: (i) $|Y|:=\sum_{i=1}^{d}Y_i$ is a $\operatorname{gamma}(|\alpha|,\beta)$
random variable, and
(ii) $X$ is independent of $|Y|$ and has Dirichlet distribution with
parameter $\alpha.$

\subsubsection{Dirichlet as a right-neutral distribution}

Let $X=(X_1,\ldots,X_d)$ be a random distribution on $\{1,\ldots,d\}$
with Dirichlet distribution $D_\alpha,\alpha\in\mathbb{R}_{+}^{d}.$
Consider the random cumulative frequencies $S_j:=\sum_{i=1}^jX_i,\
j=1,\ldots,d-1.$ Then the increments
%
\begin{equation}\label{ram}
B_j:=\frac{X_j}{1-S_{j-1}},\qquad j=1,\ldots,d-1,
\end{equation}
are independent random variables, each with a beta distribution with
parameters $(\alpha_j,|\alpha|-\sum_{i=1}^j\alpha_i)$. This
property is also known as \emph{right-neutrality} \cite{D74}.
%
Notice that such a structure holds, with different parameters, for any
reordering of the atoms of $X$.

\subsection{Size-biased Dirichlet frequencies and limit distributions}


One remarkable advantage of considering unordered versions of
Dirichlet frequencies is that they admit sensible limits as the
dimension $d$ grows to infinity, whereas the original Dirichlet
distribution is obviously bounded to finite dimensions. Two
possible ways of unordering the Dirichlet atoms are equivalent:
(1) Rearranging the frequencies in a size-biased random order; (2)
Ranking them in order of magnitude. For Dirichlet measures,
size-biased frequencies are much more mathematically treatable
than the ranked ones.

\subsubsection{Size-biased order and the GEM distribution}\label{sec:gempp}

 Let $x$ be a point of $\Delta_{(d-1)}.$ Then $x$ induces a probability
distribution on the group $\mathcal{G}_{d}$ of all permutations of $\{
1,\ldots,d\}$:
\[
\sigma_{x}(\pi)=\prod_{i=1}^{d-1}\frac{x_{\pi_{i}}}{1-\sum
_{j=1}^{i-1}x_{\pi_{j}}},\qquad \pi\in\mathcal{G}_d.
\]
Let $\alpha\in\mathbb{R}_{+}^{d}.$ The \emph{size-biased measure} on
$\Delta_{(d-1)}$ induced by a Dirichlet distribution $D_{\alpha}$ is
given by
\[
\ddot{D}_{\alpha}(A)=\int\sigma_{x}(\pi\dvtx \pi x\in A)D_{\alpha}(\mathrm{d}x).
\]
Note that $\widetilde{\sigma}_{x}\{y\}:=\sigma_{x}(\pi\dvtx \pi x=y)$ is
non-zero if and only if $y$ is a permutation of $x$, and that
\[
\widetilde{\sigma}_{x}\{y\}=\widetilde{\sigma}_{\pi x}\{y\}
=:\widetilde{\sigma}\{y\}
\qquad \forall\pi\in\mathcal{G},
\]
hence the density of the size-biased measure is
\[
\frac{\mathrm{d}\ddot{D}_{\alpha}}{\mathrm{d}y}(y)=\widetilde{\sigma}\{y\}\sum_{\pi
\in
\mathcal{G}_D}D_{\alpha}(d(\pi^{-1}y)).
\]
In particular, if $\alpha=(|\theta|/d,\ldots,|\theta|/d)$ for some
$|\theta|>0$ (symmetric Dirichlet), then its size-biased measure is
%
\begin{eqnarray}\label{sbd}
\ddot{D}_{|\theta|,d}(\mathrm{d}x)&=&d!\prod_{i=1}^{d-1}\frac{x_{i}}{1-\sum
_{j=1}^{i-1}x_{j}}D_{\alpha}(\mathrm{d}x)
\\\label{sbd2}
&\propto&\prod_{i=1}^{d-1}b_{i}^{|\theta|/d}(1-b_{i})^{
{((d-i)/d)}\theta-1}\,\mathrm{d}b_{i},
\end{eqnarray}
where $b_i=x_i/(1-\sum_{j=1}^{i-1}x_j),\ i=1,\ldots,d-1.$ So if
$\ddot{X}^{(d)}$ has distribution $\ddot{D}_{|\theta|,d}$, then
\[
\ddot{X}^{(d)}\stackrel{d}{=}\bigl(\ddot{B}_1^{(d)},\ldots,\ddot{B}_{d-1}^{(d)}\bigr),
\]
where $(\ddot{B}_{i}^{(d)})$ are $d-1$ independent beta random
variables with parameters, respectively,
$(|\theta|/d+1,({d-i}/{d})\theta),\ i=1,\ldots,d-1.$

The measure $\ddot{D}_{|\theta|,d}$ is, again, a right-neutral
measure.

Now, let $d\rightarrow\infty.$ Then $\ddot{D}_{|\theta
|,d}$ converges to the law of a right-neutral sequence
$\ddot{X}^{\infty}=(\ddot{X}_1,\ddot{X}_2,\ldots)$ such that
%
\begin{equation}\label{gem}
\ddot{X}_{j}\stackrel{\mathcal{D}}{=}\ddot{B}_{j}\prod
_{i=1}^{j-1}(1-\ddot{B}_{i}), \qquad j\geq1,
\end{equation}
for a sequence $\ddot{B}=(\ddot{B}_1,\ddot{B}_2,\ldots)$ of
independent and identically distributed ({i.i.d}.) beta weights
with parameter $(1,|\theta|)$ (here and in the following pages,
$\mathcal{D}$ means \lq\lq\emph{in distribution}\rq\rq).

\begin{definition}
The random sequence $\ddot{X}^{\infty}$ satisfying (\ref{gem}) for a
sequence of $\operatorname{beta} (1,|\theta|)$ weights is called the GEM distribution
with parameter $|\theta|$ ($\operatorname{GEM}(|\theta|)$).

Poisson point process construction \cite{K75}.
\end{definition}

 Let $Y^{\infty}=(Y_1,Y_2,\ldots)$ be the sequence
of points of a non-homogeneous point process with intensity measure
\[
N_{|\theta|}(y)= |\theta|y^{-1}\mathrm{e}^{-y}.
\]
The probability generating functional is
%
\begin{equation}\label{pdmgf}
\mathcal{F}_{|\theta|}(\xi)=\mathbb{E}_{|\theta|}\biggl(\exp
\biggl\{
\int\log\xi(y)N_{|\theta|}(\mathrm{d}y)\biggr\}\biggr)=\exp\biggl\{
|\theta|\int_{0}^{\infty}\bigl(\xi(y)-1\bigr)y^{-1}\mathrm{e}^{-y}\,\mathrm{d}y\biggr\}
\end{equation}
for suitable functions $\xi\dvtx \mathbb{R}\rightarrow[0,1].$ The
$\operatorname{GEM}(|\theta|)$ distribution can be
redefined in terms of the same point process
$Y^{\infty}$: Reorder the jumps by their size-biased
random order, that is, set
\[
\ddot{Y}_1=Y_{i1}
\]
with probability $Y_{i1}/|Y^{\infty}|$ and
\[
\mathbb{P}(\ddot{Y}_{k+1}=Y_{i,k+1}\ |
\ddot{Y}_1,\ldots,\ddot{Y}_k)=\frac{Y_{i,k+1}}{|Y|-\sum
_{j=1}^{k}\ddot{Y}_j},\qquad k=1,2,\ldots.
\]
Denote the vector of all the size-biased jumps by
$\ddot{Y}^{\infty}.$ Then
$|\ddot{Y}^{\infty}|\stackrel{\mathcal{D}}{=}|Y^{\infty}|$ is a
$\operatorname{gamma} (\theta)$ random variable,
independent of the normalized sequence
\[
\ddot{X}^{\infty}:=\frac{\ddot{Y}^{\infty}}{|\ddot{Y}^{\infty}|}
\]
and
$\ddot{X}^{\infty}$ has the $\operatorname{GEM}(|\theta|) $ distribution.

To intuitively convince oneself of such a statement, just notice that
the probability generating
functional of $\gamma_{\alpha,1}^{d},$ for
$\alpha=(|\theta|/d,\ldots,|\theta|/d),$ is \cite{G79}
%
\begin{eqnarray}\label{pdlim}
\mathcal{F}_{|\theta|,d}(\xi)&=&\biggl(\int_0^{\infty}\xi
(y)\gamma_{{|\theta|/d},1}(\mathrm{d}y)\biggr)^{d}\nonumber
\\
&=&\biggl(1+\int_0^\infty\bigl(\xi(y)-1\bigr)\frac{|\theta|}{d}\frac
{y^{{|\theta|/d}-1}\mathrm{e}^{-y}}{\Gamma({|\theta
|/d}+1)}\,\mathrm{d}y\biggr)^{d}
\\
&\mathop{\rightarrow}\limits_{d\rightarrow\infty}&
\mathcal{F}_{|\theta|}(\xi) ,\nonumber
\end{eqnarray}
so a finite size-biased collection of $d$ {i.i.d.}, normalized
gamma jumps has a $\operatorname{GEM}(\theta)$ limit distribution, as $d\rightarrow
\infty$.

\subsubsection{Beta-Stacy distributions}

The measures $D_{\alpha},\ddot{D}_{|\theta|,d},\ddot{D}_{|\theta
|}$ are all right-neutral distributions with independent beta parameters.

\begin{definition}
For $d\leq\infty,$ let $B^*_1,\ldots,B^*_{d-1}$ be
a collection of mutually independent beta random variables with
parameters $\{\alpha_i,\beta_i\}_{i=1}^{d}$ (if $d=\infty$, take an
infinite sequence of such weights). A random discrete distribution
$X\in\Delta_{(d-1)}$ is said to have a \emph{beta-Stacy} law if
$X_1\stackrel{\mathcal{D}}{=}B^*_1$ and, for every $j\leq d-1$,
%
\[
1-\sum_{i=1}^{j-1}X_i\stackrel{\mathcal{D}}{=}\prod_{i=1}^{j-1}(1-B_i^*).
\]
A notable example of infinite-dimensional beta-Stacy distribution is
the \emph{two-parameter $\operatorname{GEM}(\alpha, \theta)$ distribution \cite
{P95,P96}} whereby, for every $j\leq d-1,$ $B^*_j$ is a $\operatorname{beta}(1-\sigma
,\theta+j\sigma)$ random variable, with either $\sigma\in[0,1]$ and
$\theta>-\sigma$ or $\sigma<0$ and $\theta=|\sigma|m$ for some
$m\in\mathbb{N}$.

The two-parameter GEM distribution is the most general class of
right-neutral distributions that is also invariant under size-biased
permutation; other remarkable properties (it is regenerative and \emph
{Gibbs}) make it one of the most studied models for generating
consistent, exchangeable random partitions (see \cite{P06} and
references therein).
\end{definition}

\subsection{Sampling formulae}

The multinomial-Dirichlet distribution can be obtained by mixing the
parameter of a~multinomial distribution with a Dirichlet mixing
measure: If $X$ has $D_\alpha$ distribution,
%
\begin{eqnarray}\label{dm2}
\mathit{DM}_{\alpha}(r;|r|)=\mathbb{E}\biggl[\pmatrix{{|r|}\cr r }X^r
\biggr]
=\pmatrix{{|r|}\cr  r}\frac{\prod_{i=1}^{d}(\alpha_i)_{(r_i)}}{(|\alpha
|)_{(|r|)}},
\end{eqnarray}
where $(a)_{(x)}:=\Gamma(a+x)/\Gamma(a)$ for $a>0$.

\subsubsection{Partial right-neutrality}

For every $r\in\mathbb{N}^{d}$ and $\alpha\in\mathbb{R}^{d}_+$, denote
as usual $R_j=\sum_{i=j+1}^{d}r_i$ and
$A_j=\sum_{i=j+1}^{d}\alpha_i$. It is easy to see that
\begin{eqnarray}\label{mdneut}
\mathit{DM}_{\alpha}(r;R)&=&\prod_{j=1}^{d-1}\pmatrix{R_{j-1}\cr r_{j}}\int_0^1
z_j^{r_j}(1-z_j)^{R_j}D_{\alpha_j,A_j}(\mathrm{d}z_j)\nonumber
\\[-8pt]\\[-8pt]
&=&\prod_{j=1}^{d-1}\mathit{DM}_{\alpha_j,A_j}(r_j;R_{j-1}).\nonumber
\end{eqnarray}
In other words, for every $j=1,\ldots,d-1$, ${r_j}/{R_j}$ is \emph
{conditionally} independent of $r_1,\ldots,r_{j-1}$, given $R_j$. Such a
property, a direct consequence of the
Dirichlet, is responsible for our construction of multivariate Hahn polynomials.

\subsubsection{Negative binomial sums}\label{snbsum}

Another construction of $\mathit{DM}_{\alpha}$ is possible,
based on negative binomial random sequences, which parallels the
gamma construction of the
Dirichlet measure of Section \ref{sec:gammasums}.

Let $\mathit{NB}_{|\alpha|,y}(k)\dvtx |\alpha|>0,$ denote the \emph{negative
binomial distribution} with probability mass function:
%
\begin{equation}
\mathit{NB}_{|\alpha|,p}(k)=\frac{{(|\alpha
|)}_{(k)}}{k!}p^{k}(1-p)^{|\alpha|},\qquad k=0,1,\ldots.
\label{nbin}
\end{equation}
With both parameters in $\mathbb{N}$, such a measure describes the
distribution of the number of failures occurring in a sequence of
{i.i.d.} Bernoulli experiments (with success probability $1-p$),
before the $\alpha$th success.

Two features of $\mathit{NB}_{|\alpha|,p}$
will prove useful, in Section
\ref{sec:meix}, to connect multiple Meixner polynomials to
multivariate Hahn polynomials.

\begin{enumerate}[(1)]
\item[(1)] Poisson--gamma mixtures:
\begin{eqnarray}\label{pgd}
\mathit{NB}_{|\alpha|,p}(k)&=&\int_0^\infty Po_{\lambda}(k)\gamma_{|\alpha
|,{p/(1-p)}}(\mathrm{d}\lambda),\nonumber
\\[-8pt]\\[-8pt]
\mathit{Po}_{\lambda}(k)&=&\frac{\lambda^{k}\mathrm{e}^{-\lambda}}{k!},\qquad
k=0,1,2,\ldots.\nonumber
\end{eqnarray}

\item[(2)] Normalized negative binomial vectors.
\end{enumerate}

Consider any $\alpha\in\mathbb{R}_{+}^{d}$ and $p\in(0,1).$ Let
$R_1,\ldots,R_d$ be independent negative binomial random variables
with parameter $(\alpha_i,p)$, respectively, for $i=1,\ldots,d$.
Then
\begin{longlist}[(ii)]
\item[(i)] $|R|:=\sum_{i=1}^{d}R_i$ has law $\mathit{NB}_{|\alpha|,p}$.

\item[(ii)] Conditional on $|R|=|r|,$ the vector $R=(R_1,\ldots,R_d)$ has a
Dirichlet-multinomial distribution with parameter $(\alpha,|r|)$:
%
\begin{equation}\label{nbsum}
\prod_{i=1}^{d}\mathit{NB}_{\alpha_i,p}(r_i)=\mathit{NB}_{|\alpha|,p}(|r|)\ \mathit{DM}_{\alpha
}(r;|r|).
\end{equation}
\end{longlist}

\subsubsection{Hypergeometric distribution}\label{sec:hg}
Consider the form of the probability mass function
$\mathit{DM}_{\alpha}$ but now replace the parameter~$\alpha$ with
$-\varepsilon=(-\varepsilon_1,\ldots,-\varepsilon_d)$ with $0\leq n_j\leq
\varepsilon_{j},\ j=1,\ldots,d.$ Then
%
\begin{eqnarray}\label{H}
\mathit{DM}_{-\varepsilon}(n)=\frac{|n|!}{n_1 !\cdots n_d!}\frac{{(-\varepsilon
)}_{(n)}}{{(-|\varepsilon|)}_{(|n|)}}
=\frac{\prod_{i=1}^{d}\matrix{{\varepsilon_i}\choose {n_i}}}{\matrix{{|\varepsilon
|}\choose {|n|}}}=:H_{\varepsilon}(n).
\end{eqnarray}
$H_{\varepsilon}(n)$ is known as the \emph{multivariate
hypergeometric distribution} with parameter $\varepsilon.$

The partial right-neutrality property of the Dirichlet-multinomial
distribution is preserved for the hypergeometric law; however, the
interpretation as a Dirichlet mixture of {i.i.d.} laws is lost
as the Dirichlet (as well as the gamma and the beta) integral is
not defined for negative parameters.

\subsection{Conjugacy properties}\label{sec:conjd}

The gamma and the Dirichlet distribution, and,
similarly, the negative binomial and the Dirichlet-multinomial
distributions, are entangled by yet
another property known in Bayesian statistics as \emph{conjugacy} with
respect to sampling.

A statistical model can be described by a probability triplet
$\{M,\mathcal{M},l_{\Lambda}\}_{\Lambda\in E}$, where the
likelihood function $l_{\Lambda}(x)$ depends on a random parameter
$\Lambda$ living in some probability space $(E,\mathcal{E},\pi).$ The
distribution $\pi$ of $\Lambda$ is called the \emph{prior} measure of
the model. The \emph{posterior} measure of the model is any
version $\pi_x(\cdot)=\pi(\cdot|X=x) $ of the conditional
probability satisfying
%
\begin{equation}\label{conj}
 \int_A \pi(B | X=x)
\int l_\lambda(\mathrm{d}x)\pi(\mathrm{d}x)=\int_B
l_\lambda(A)\pi(\mathrm{d}\lambda)\qquad\mbox{a.s. }\forall
A\in\mathcal{M}, B\in\mathcal{E}.
\end{equation}

\begin{definition}
\label{def:conj} Let $\mathcal{C}$ be a family of prior measures
for a statistical model with likelihood~$l_{\Lambda}$.
$\mathcal{C}$ is \emph{conjugate} with respect to $l_{\Lambda}$ if
\[
\pi\in\mathcal{C}\quad\Longrightarrow\quad\pi_x\in\mathcal{C} \qquad\forall x.
\]
\end{definition}

It is easy to check that both gamma and Dirichlet measures are
conjugate classes of prior measures. Bayes' theorem shows us the
role as \emph{marginal distributions} played, respectively, by
$\mathit{NB}_{\alpha,p}$ and $\mathit{DM}_{\alpha}$.

\begin{example}
The class of gamma priors is conjugate with respect to
$l_{\lambda}=\mathit{Po}_{\lambda}$ on $\{0,1,2,\ldots\}.$ The posterior
measure is
%
\begin{equation}
\pi_{x}(\mathrm{d}\lambda)=\frac{\mathit{Po}_{\lambda}(x)\gamma_{\alpha,\beta
}(\mathrm{d}\lambda)}{\mathit{NB}_{\alpha,{\beta/(1+\beta)}}(x)}
=\gamma_{\alpha+x,{\beta/(1+\beta)}}(\mathrm{d}\lambda).
\end{equation}
Similarly, the class of multivariate gamma priors $\{\gamma^d_{\alpha
,\beta}:\alpha\in\mathbb{R}^{d},\beta>0\}$ is conjugate with
respect to
$\{\mathit{Po}_{\lambda}^{d}(x),\lambda\in\mathbb{R}_{+}^{d},x\in\mathbb
{N}^{d}\}$.
\end{example}

\begin{example}
The class of beta priors
$\{D_{\alpha,\beta}\dvtx (\alpha,\beta)\in\mathbb{R}_{+}^{2}\}$ is
conjugate with respect to the binomial likelihood
$l_{\lambda}=B_{\lambda}(\cdot)$ on $\{0,1,2,\ldots,|n|\},$ for
any integer $|n|.$ The posterior distribution is
%
\begin{equation}
\pi_{r}(\mathrm{d}\lambda)=\frac{B_{\lambda}(|r|,|n-r|)D_{\alpha,\beta
}(\mathrm{d}\lambda)}{\mathit{DM}_{\alpha,\beta}
(|r|;|n|-|r|)}=D_{\alpha+|r|,\beta+|n|-|r|}(\mathrm{d}\lambda).
\end{equation}
Similarly, the class of Dirichlet measures is conjugate with
respect to multinomial sampling. 
\end{example}

%

\section{Jacobi polynomials on the simplex}\label{sec3}

If $X,Y$ are independent random variables, their distribution
$W_{X,Y}$ is the product $W_{X}W_{Y}$ of their marginal
distributions, and therefore orthogonal polynomials $Q_{n,k}(x,y)$
in $W_{X,Y}$ are simply obtained by products $P_{n}(x)R_{k}(y)$ of
orthogonal
polynomials with $W_X$ and~$W_Y$ as weight measures, respectively.

The key idea for deriving multivariate polynomials with respect to
Dirichlet measures on the simplex, and to all related
distributions treated in the subsequent sections, exploits the
several properties of \emph{conditional} independence enjoyed by
the increments of~$D_{\alpha}$, as pointed out in Section
\ref{sec:gammasums}. A method for constructing orthogonal
polynomials in the presence of a particular kind of conditional
independence, where $Y$ depends on $X$ only through a polynomial
$\rho(x)$ of the first-order, is illustrated by the following multidimensional
modification of Koornwinder's method (see \cite{Koo75}, Section~3.7.2).

\begin{proposition}
\label{prop:mvkoo} For $l,d\in\mathbb{N},$ let $(X,Y)$ be a random
point of $\mathbb{R}^{l}\times\mathbb{R}^{d}$ with distribution~$W$. Let
$\rho:\mathbb{R}^{l}\rightarrow\mathbb{R}$ define
polynomials on $\mathbb{R}^{l}$ of order at most 1.

Assume that the random variable
\[
Z:=\frac{Y}{\rho(X)}
\]
is independent of $X.$ Denote with $W_X$ and $W_Z$ the marginal
distributions of $X$ and $Z$, respectively. Then a system of multivariate
polynomials, orthogonal with respect to $W,$ is given by
\begin{eqnarray} \label{mvkoo}
&&G_{n}(x,y)=P^{(N_{l})}_{(n_1,\ldots,n_l)}(x)(\rho
(x))^{N_l}R_{(n_{l+1},\ldots,n_{l+d})}\biggl(\frac{y}{\rho(x)}
\biggr),\nonumber
\\[-8pt]\\[-8pt]
&&\quad(x,y)\in\mathbb{R}^{l}\times\mathbb{R}^{d},
n\in\mathbb{N}^{l+d},\nonumber
\end{eqnarray}
where $N_{l}=n_{l+1}+\cdots+ n_{l+d}$, $\{P_{k}^{(|m|)}\}_{k\in
\mathbb{R}^{l}}$ and $\{R_{m}\}_{m\in\mathbb{R}^{d}}$ are systems of
orthogonal
polynomials with weight measures given by $(\rho(x))^{2|m|}W_X$ and
$W_Z,$ respectively.
\end{proposition}

\begin{pf}
When $d=l=1$ this proposition is essentially a probabilistic
reformulation of Koornwinder's construction (\cite{Koo75}, Section~3.7.2).
The proof is similar for any $l,d$. That~$G_{n}$ is a polynomial
of degree $|n|$ is evident as the denominator of the term of
maximum degree in $R$ simplifies with
$(\rho(x))^{n_{l+1}+\cdots+n_{l+d}}.$ To show orthogonality, note
that the assumption of conditional independence implies that
\[
W(\mathrm{d}x,\mathrm{d}y)=W_X(\mathrm{d}x)W_Z\biggl(\frac{1}{(\rho(x))^{d}}\,\mathrm{d}y\biggr).
\]
Denote $b_n=\mathbb{E}[P_n^2]$ and $c_n=\mathbb{E}[R_n^2]$,
$n=0,1,2,\ldots.$ For $k,r\in\mathbb{R}^{l}$ and $m,s\in\mathbb{R}^{d},$
\begin{eqnarray*}
&&\int G_{(k,m)}(x,y)G_{(r,s)}(x,y)W(\mathrm{d}x,\mathrm{d}y)
\\
&&\quad=\int
P^{m}_{k}(x)P^{s}_{r}(x)(\rho(x))^{m+s}W_X(\mathrm{d}x)\int
R_{m}(z)R_{s}(z)W_Z(\mathrm{d}z)\nonumber
\\
&&\quad=\int P^{m}_{k}(x)P^{m}_{r}(x)(\rho(x))^{2m}W_X(\mathrm{d}x)c_{m}\delta
_{ms}
\\
&&\quad=b_{k}c_m\delta_{kr}\delta_{ms}.\nonumber
\end{eqnarray*}
\upqed\end{pf}

\subsection{$d=2$; Jacobi polynomials on $[0,1]$}
\label{*}

For $d=2$, $D_\alpha$ reduces to the beta
distribution, the weight measure of (shifted) Jacobi polynomials.
These are functions of one variable living in
$\Delta_1\equiv[0,1]$. It is convenient to recall some known
properties of such polynomials. Consider the measure
%
\begin{equation}\label{lambda}
\widetilde{w}_{a,b}(\mathrm{d}x)=(1-x)^{a}(1+x)^{b}\mathbb{I}\bigl(x\in
(-1,1)\bigr)\,\mathrm{d}x,\qquad a,b>-1,
\end{equation}
where $\mathbb{I}(A)$ is the indicator function, equal to 1 if
$A$, and $0$ otherwise. This is the weight measure of the Jacobi
polynomials defined by
\[
\widetilde{P}_{n}^{a,b}(x):=\frac{{(a+1)}_{(n)}}{n!}{}_{2}F_1\left(
\begin{array}{cc}-n,\ n+a+b+1\\
a+1
\end{array}
\ \vline\ \ \frac{1-x}{2}\right),
\]
where
${}_{p}F_q, p,q\in\mathbb{N},$ denote the hypergeometric
function (see \cite{AS} for basic properties).

The normalization constants are given by the relation
%
\begin{equation}\label{nc}
\int_{(-1,1)}\widetilde{P}_{n}^{a,b}(x)\widetilde
{P}_{m}^{a,b}(x)\widetilde{w}_{a,b}(\mathrm{d}x)=\frac{2^{a+b+1}}{2n+a+b+1}
\frac{\Gamma(n+a+1)\Gamma(n+b+1)}{n!\Gamma(n+a+b+1)}\delta
_{mn}.
\end{equation}

The Jacobi polynomials are known to be solution of the second-order
partial differential equation
%
\begin{equation}\label{jpde}
(1-x^2)y^{\prime\prime}(x)+[b-a-x(a+b+2)]y^{\prime}(x)=-n(n+a+b+1)y(x).
\end{equation}
By a simple shift of measure it is easy to see that, for
$\alpha,\beta>0$ and $\theta:=\alpha+\beta$, the modified
polynomials
%
\begin{equation}\label{01J}
P_{n}^{\alpha,\beta}(x)=\frac{n!}{{(n+\theta-1)}_{(n)}}\widetilde
{P}_{n}^{\beta-1,\alpha-1}(2x-1),\qquad{\alpha,\beta>0},
\end{equation}
are orthogonal with respect to
the beta distribution on $[0,1],$ which can be written as
%
\begin{equation}\label{sbeta}
D_{\alpha,\beta}(\mathrm{d}x)=\frac{\widetilde{w}_{\beta
-1,\alpha
-1}(\mathrm{d}u)}{2^{\alpha+\beta-1}B(\alpha,\beta)},
\end{equation}
where $u=2x-1$.

 Denote the standardized Jacobi polynomials with
\[
\widetilde{R}_{n}^{a,b}(x)=\frac{\widetilde
{P}_{n}^{a,b}(x)}{\widetilde{P}_{n}^{a,b}(1)}
\quad\mbox{and} \quad
R_{n}^{\alpha,\beta}(x)=\frac{P_{n}^{\alpha,\beta
}(x)}{P_{n}^{\alpha,\beta}(1)}.
\]
Obviously
%
\begin{equation}\label{01nJ}
R_{n}^{\alpha,\beta}(x)=\widetilde{R}_{n}^{(\beta-1,\alpha-1)}(2x-1).
\end{equation}
By (\ref{nc}) the new constant of proportionality is
%
\begin{eqnarray}\label{hn}
\frac{1}{\zeta_{n}^{(\alpha,\beta)}}&:=&\int_{0}^{1}[R_{n}^{\alpha
,\beta}(x)]^{2}D_{\alpha,\beta}(\mathrm{d}x)\nonumber
\\[2pt]
&=&\biggl(\frac{{(\theta+n-1)}_{(n)}}{{(\beta)}_{(n)}}
\biggr)^{2}\frac{n!
{\alpha}_{(n)}{(\beta)}_{(n)}}{{(\theta)}_{(2n)}{(\theta+n-1)}_{(n)}}
\\[2pt]
&=&n!\frac{1}{(\theta+2n-1){(\theta)}_{(n-1)}}\frac{
{(\alpha)}_{(n)}}{{(\beta)}_{(n)}}, \qquad n=0,1,\ldots.\nonumber
\end{eqnarray}
A symmetry relation is
%
\begin{equation}\label{symj}
{R}_{n}^{\alpha,\beta}(x)=\frac{{R}_{n}^{\beta,\alpha
}(1-x)}{{R}_{n}^{\beta,\alpha}(0)}.
\end{equation}
Note that, if $\{{P_{n}^{*}}^{\alpha,\beta}(x)\}$ is a
system of \emph{orthonormal} polynomials with weight measure~%
$D_{\alpha,\beta}$, then
%
\begin{equation}
\zeta_{n}^{(\alpha,\beta)}=[{P_{n}^{*}}^{\alpha,\beta}(1)]^{2}.
\end{equation}

\subsection{$2\leq d<\infty$. Multivariate Jacobi polynomials on the simplex from right-neutrality}

 A system of multivariate polynomials with respect to a
Dirichlet distribution on $d\leq\infty$ points can be derived by
using its right-neutrality property, via Proposition
\ref{prop:mvkoo}. Let
$\mathbb{N}_{d,|m|}=\{{n}=(n_{1},\ldots,n_{d})\in
\mathbb{N}^{d}\dvtx |{n}|=|m|\}$. For every $n\in\mathbb{N}_{d-1,|n|}$
and $\alpha\in\mathbb{R}_{+}^{d}$ denote
$
N_{j}=\sum_{i=j+1}^{d-1}n_{i}$ and $A_{j}=\sum_{i=j+1}^{d}\alpha_{i}.$

\begin{proposition}\label{mop}
For $d<\infty$, a system of multivariate orthogonal polynomials
on the Dirichlet distribution $D_{\alpha}$ is given by
%
\begin{equation}\label{gJ}
R_{n}^{\alpha}(x)=\prod_{j=1}^{d-1}R_{n_{j}}^{(\alpha_{j},A_{j}+2N_{j})}
\biggl(\frac{x_{j}}{1-s_{j-1}}\biggr)(1-s_{j-1})^{N_{j}},\qquad x\in\Delta_{(d-1),}
\end{equation}
where $s_{j}=\sum_{i=1}^{j}x_{i}$.
\end{proposition}

 Notice that
$
R_{n}^{\alpha}(\mathbf{e}_{d})=1, $
where $\mathbf{e}_{j}:=(\delta_{ij}\dvtx i=1,\ldots,d)$.
A similar definition for polynomials in the Dirichlet
distribution is proposed in \cite{KS91}, in terms of non-shifted
Jacobi polynomials
$\widetilde{R}_{n}$. For an alternative choice of basis, see \cite
{DX02}.\vadjust{\goodbreak}

\begin{pf*}{Proof of Proposition \ref{mop}}
The polynomials in $R_{n}^{\alpha}(x)$ given in Proposition \ref{mop}
admit a recursive definition as follows:
\begin{eqnarray}\label{recdJ}
&&R_{n_1,\ldots,n_{d-1}}^{\alpha}(x_{1},\ldots,x_{d})\nonumber
\\[-8pt]\\[-8pt]
&&\quad=
R_{n_1}^{(\alpha
_{1},A_{1}+2N_{1})}(x_{1})(1-x_{1})^{N_{1}}R_{n_2,\ldots
,n_{d-1}}^{\alpha^{*}_{2}}
\biggl(\frac{x_{2}}{1-x_{1}},\ldots,\frac{x_{d}}{1-x_{1}}\biggr),\nonumber
\end{eqnarray}
where $\alpha^{*}_{j}=(\alpha_j,\ldots,\alpha_ {d})$ ($j\leq
d-1$); so Proposition \ref{prop:mvkoo} is used with $l=1,
\rho(x)=1-x$ and inductively on $d$. The claim is a consequence of
the neutral-to-the-right property and Proposition~\ref{prop:mvkoo} -- for consider the orthogonality of a term
%
\begin{equation}\label{jterm}
\biggl(1-\frac{X_{j}}{1-S_{j-1}}\biggr)^{N_{j}}R_{n_{j}}^{\alpha
_{j},A_{j}+2N_j}\biggl(\frac{X_{j}}{1-S_{j-1}}\biggr)
\end{equation}
in
$R_{n}^{\alpha}$ with a similar term in $R_{m}^{\alpha}$ for some
$m=(m_{1},\ldots,m_{d-1})$-polynomial. Assume without loss of
generality that
for some $j=1,\ldots,d-1$, $m_{k}=n_{k}$ for $k=j+1,\ldots,d-1$ and
$m_{j}< n_{j}$. Then $N_{j}=M_{j}$ and, multiplying the product of
(\ref{jterm}) by the corresponding beta density $D_{\alpha
_j,A_j}(\mathrm{d}B_j)/\mathrm{d}B_j,$ where $B_{j}$ is as in (\ref{ram}), gives
%
\begin{equation}\label{t1}
B_{j}^{\alpha_{j}-1}(1-B_{j})^{A_{j}+2N_{j}-1}R_{n_{j}}^{\alpha
_{j},A_{j}+2N_{j}}(B_{j}) R_{m_{j}}^{\alpha_{j},A_{j}+2N_{j}}(B_{j}).
\end{equation}
 Since $R_{n_{j}}$ is orthogonal to polynomials of degree
less than $n_{j}$ on the weight measure
$D_{\alpha_{j},A_{j}+2N_{j}}$, then the integral with respect to
$\mathrm{d}B_{j}$ of the quantity (\ref{t1}) vanishes, which proves the
orthogonality.
\end{pf*}

 The orthogonality constant for $\{R^{\alpha}_n\}$ can be
easily derived as
\begin{eqnarray}\label{mopconst}
\frac{1}{\zeta^{\alpha}_{n}}&:=&\int_{\Delta_{(d-1)}}
(R^{\alpha}_n(x))^{2}D_{\alpha}(\mathrm{d}x)=\frac{1}{\prod
_{j=1}^{d-1}\zeta^{\alpha_{j},A_{j}+2N_{j}}_{n_{j}}}\nonumber
\\[-8pt]\\[-8pt]
&=&\prod_{j=1}^{d-1}\frac{n_j!
{(\alpha
_j)}_{(n_j)}}{{(A_{j-1}+N_{j})}_{(n_j-1)}(A_{j-1}+2N_{j-1}-1){(A_j+2N_j)}_{(n_j)}}.\nonumber
\end{eqnarray}
 Notice that the same construction shown in Proposition \ref
{mop} could be similarly expressed in terms of the polynomials
$\{P_{n_{j}}^{\alpha_{j},A_{j}+2N_{j}}\}$ or
$\{{P^{\star}}^{\alpha_{j},A_{j}+2N_{j}}\}$ instead of
$\{R_{n_{j}}^{\alpha_{j},A_{j}+2N_{j}}\},$ the only difference
resulting in the orthogonality constants.

\subsection{Multivariate Jacobi on beta-Stacy distributions} \label{sec:rank}

Random distributions of beta-Stacy type are all
right-neutral. Orthogonal polynomials with respect to general
beta-Stacy measures can be therefore constructed in very much the same
way as in Proposition \ref{mop}, with a similar proof.

\begin{proposition}
Let $d\leq\infty$ and $(\alpha,\beta)\in\mathbb{R}_+^d\times
\mathbb{R}_+^d. $ Let $\mu_{\alpha,\beta}$ be the distribution of a~%
$\operatorname{beta\mbox{-}Stacy}(\alpha,\beta)$ random point of $\Delta_{(d-1)}.$ A
system of orthogonal polynomials in $\mu_{\alpha,\beta}$ is given by
%
\begin{equation}\label{bspoly}
{R}_{n}^{*(\alpha,\beta)}(x)=\prod
_{j=1}^{d-1}R_{n_{j}}^{(\alpha_j,\beta_j+2N_{j})}
\biggl(\frac{x_{j}}{1-s_{j-1}}\biggr)(1-s_{j-1})^{N_{j}},\qquad
x\in\Delta_{(d-1)},n\in\mathbb{N}^{d}.
\end{equation}
The constant of orthogonality is given by
\begin{eqnarray}\label{bsnc}
\frac{1}{\zeta^{\alpha,\beta}_{n}}&=&\frac{1}{\prod
_{j=1}^{d-1}\zeta^{\alpha_{j},\beta_{j}+2N_{j}}_{n_{j}}}\nonumber
\\[-8pt]\\[-8pt]
&=&\prod_{i=1}^{d-1}\frac{n_i!{(\alpha_i)}_{(n_i)}}{(\alpha
_i+\beta_i+2N_{i-1}-1){(\alpha_i+\beta_i+2N_i)}_{(n_i-1)}{(\beta
_i+2N_i)}_{(n_i)}}.\nonumber
\end{eqnarray}
\end{proposition}

\begin{example}
We have seen that all size-biased Dirichlet measures are beta-Stacy. A~%
system of orthogonal polynomials in $\ddot{D}_{|\theta|,d}$ is
\begin{eqnarray}\label{dsbdir}
&&\ddot{R}_{n}^{(|\theta|,d)}(x)=\prod
_{j=1}^{d-1}R_{n_{j}}^{(|\theta|/d+1,{((d-j)/d)}\theta+2N_{j})}
 \biggl(\frac{x_{j}}{1-s_{j-1}}\biggr)(1-s_{j-1})^{N_{j}},\nonumber
\\[-8pt]\\[-8pt]
 &&\quad x\in\Delta_{(d-1)},n\in\mathbb{N}^{d}.\nonumber
\end{eqnarray}
\end{example}

\begin{example}
As $d\rightarrow\infty,$ $\ddot{D}_{|\theta|,d}$
converges to the so-called $\operatorname{GEM}(\theta)$ distribution, that is, an
infinite-dimensional beta-Stacy with all i.i.d.
weights being beta random variables with parameter $(\alpha_j,\beta
_j)=(1,\theta).$
Let
$\ddot{D}_{|\theta|,\infty}=\lim_{d\rightarrow\infty}\ddot
{D}_{|\theta|,d}$
denote the GEM distribution with parameter $|\theta|.$ For $|\theta
|>0$, an orthogonal system with
respect to the weight measu\-re~$\ddot{D}_{|\theta|,\infty}$ is
given by the polynomials:
\begin{eqnarray}\label{gemp}
&&\ddot{R}_{n}^{|\theta|}(x)=\prod_{j=1}^{\infty
}R_{n_{j}}^{(1,\theta+2N_{j})}
\biggl(\frac{x_{j}}{1-s_{j-1}}\biggr)(1-s_{j-1})^{N_{j}},\nonumber
\\[-8pt]\\[-8pt]
&&\quad x\in\Delta_{\infty},n\in\mathbb{N}^{\infty}\dvtx |n|=0,1,\ldots.\nonumber
\end{eqnarray}
\end{example}

\begin{example}
For the two-parameter $\operatorname{GEM}(\sigma,\theta)$ distribution, $\alpha
_j=1-\sigma$ and $\beta_j=\theta+j\sigma$. The polynomials are of
the form
\begin{eqnarray}\label{2gemp}
&&\ddot{R}_{n}^{\sigma,\theta}(x)=\prod_{j=1}^{\infty
}R_{n_{j}}^{(1-\sigma,\theta+j\sigma+2N_{j})}
\biggl(\frac{x_{j}}{1-s_{j-1}}\biggr)(1-s_{j-1})^{N_{j}},\nonumber
\\[-8pt]\\[-8pt]
&&\quad x\in\Delta_{\infty},n\in\mathbb{N}^{\infty}\dvtx |n|=0,1,\ldots.\nonumber
\end{eqnarray}
\end{example}

\section{Multivariate Jacobi and multiple Laguerre polynomials} \label{sec:lag}
The Laguerre polynomials, defined by
%
\begin{equation}\label{lag}
L_{|n|}^{|\alpha|}(y)=\frac{{(|\alpha|)}_{(|n|)}}{|n|!}
{}_{1}F_{1}(-|n|;|\alpha|;y),\qquad |\alpha|>0,
\end{equation}
are orthogonal to the gamma density $\gamma_{|\alpha|,1}$ with
constant of orthogonality
%
\begin{equation}\label{lcon}
\int_0^{\infty}\bigl[L_{|n|}^{|\alpha|}(y)\bigr]^{2}\gamma_{|\alpha
|}(\mathrm{d}y)=\frac{{(|\alpha|)}_{(|n|)}}{|n|!}.
\end{equation}
(Note that the usual convention is to define Laguerre polynomials in
terms of the parameter $|\alpha'|:=|\alpha|-1>-1.$ Here we prefer to use
positive parameter for consistency with the parameters in
the amma distribution.)

\begin{remark}\label{rmscale}
If $Y$ is a $\operatorname{gamma} (|\alpha|)$ random
variable, then, for every
scale parameter $\beta\in\mathbb{R}_{+}$, the distribution of
$Z:=\beta
Y$ is $\gamma_{|\alpha|,\beta}(\mathrm{d}z).$ Thus the system
\[
\biggl\{L_n^{|\alpha|}\biggl(\frac{z}{\beta}\biggr)\biggr\}
_{n=0,1,\ldots}
\]
is orthogonal with weight measure $\gamma_{|\alpha|,\beta}.$
\end{remark}

Let $Y\in\mathbb{R}_{+}^{d}$ be a random vector with distribution
$\gamma_{\alpha,\beta}^{d}.$
By the stochastic independence of its coordinates, orthogonal
polynomials of degree $|n|$ with the distribution of $Y$ as weight
measure are
simply
%
\begin{equation}\label{mlag1}
L^{\alpha,
\beta}_{n}(y)=\prod_{i=1}^{d}L_{n_{i}}^{\alpha_i}\biggl(\frac
{y_{i}}{\beta}\biggr),\qquad y\in\mathbb{R}^{d},n\in\mathbb
{N}_{n},
\end{equation}
with constants of orthogonality of
%
\begin{equation}\label{mlagcp}
\frac{1}{\varphi_n}=\mathbb{E}(L^{\alpha}_{n}(Y)
)^{2}=\prod_{i=1}^{d}\frac{{(\alpha_i)}_{(n_i)}}{n_i!}.
\end{equation}
Therefore, with the notation introduced in Section \ref
{sec:gammasums}, because of the one-to-one mapping
\[
(Y_1,\ldots,Y_d)\mapsto(|Y|,X_1,\ldots,X_d),
\]
one can obtain an alternative system of orthogonal polynomials from
$y_1,\ldots,y_n.$

\begin{proposition}
\label{pr:mlag}
The polynomials defined by
%
\begin{equation}\label{mlag}
L^{\alpha, \beta*}
_{n}(y)=L_{n_d}^{|\alpha|+2|n'|}\biggl(\frac{|y|}{\beta}
\biggr)\biggl(\frac{|y|}{\beta}\biggr)^{|n'|}R^{\alpha}_{n'}\biggl(\frac
{y}{|y|}\biggr),\qquad n\in\mathbb{N}^{d},y\in\mathbb{R}^{d},
\end{equation}
with $n'=(n_1,\ldots,n_{d-1})$ and $R^{\alpha}_{m}$ defined by (\ref
{gJ}), are orthogonal with respect to $\gamma_{\alpha, \beta}^{d}.$
\end{proposition}

\begin{pf}
The proof of (\ref{mlag}) is straightforward and follows
immediately from Proposition~\ref{prop:mvkoo}, with $l=1$, $X=|Y|$ and
$\rho(x)=x$ (remember that $|Y|$ is gamma with parameter $(|\alpha
|,\beta)$).
\end{pf}

 From now on we will only consider the case with
$\beta=1$, without much loss of generality. The constant of
orthogonality of the resulting system $\{L^{\alpha*}_n\}$ is
%
\begin{eqnarray}\label{mlconst}
\frac{1}{\varphi_{n}^{*}}&:=&\int_{\mathbb{R}^{d}}[L^{\alpha
*}_{n}(y)]^2\prod_{i=1}^{d}\gamma_{\alpha_i}(\mathrm{d}y_i)\nonumber
\\
&=&\int_{0}^{\infty}\bigl[L_{n_d}^{|\alpha
|+2(|n|-n_d)}(|y|)|y|^{|n|-n_d}\bigr]^2\gamma_{|\alpha|}(\mathrm{d}|y|)\int
_{\Delta_{(d-1)}}
[R^{\alpha}_{n'}(x)]^2 D_{\alpha}(\mathrm{d}x)\nonumber
\\
&=&\frac{{(|\alpha|)}_{(2|n'|)}}{\zeta^{\alpha}_{n'}}\int
\bigl[L_{n_d}^{|\alpha|+2|n'|}(|y|)\bigr]^2\gamma_{\alpha
+2|n'|}(\mathrm{d}|y|)\nonumber
\\
&=&\frac{1}{n_d!}\frac{({(|\alpha|)}_{(2|n'|)})^{2}}{\zeta
^{\alpha}_{n'}},
\end{eqnarray}
where $\zeta^{\alpha}_{n'}$ is as in (\ref{mopconst}).

\subsection{Connection coefficients}

 The two systems $L^{\alpha}_{n}$ and $L^{\alpha
*}_{n}$ can be expressed as linear combinations of each other:
%
\begin{equation}\label{lc1}
L^{\alpha*}_{n}(y)=\sum_{|m|=|n|}\varphi_m c^{*}_m(n) L_{m}^{\alpha
}(y)
\end{equation}
and
%
\begin{equation}\label{lc2}
L^{\alpha}_{n}(y)=\sum_{|m|=|n|}\varphi^{*}_m c_m(n) L_{m}^{\alpha
*}(y),
\end{equation}
where
\[
c^{*}_m(n)\delta_{|m||n|}=\mathbb{E}[L^{\alpha
*}_{n}(y)L^{\alpha}_{m}(y)]=c_{n}(m)\delta_{|m||n|}.
\]
 For general $m,n$ a representation for $c^{*}_m(n)$ can
be derived in terms of a mixture of Lauricella functions of the
first (A) type. Such functions are defined \cite{L1893} as
\[
F_{A}(|a|;b;c;z)=\sum_{m\in\mathbb{N}^d}\frac{1}{m_1!\cdots
m_d!}\frac
{{|a|}_{(|m|)}b_{(m)}}{{c}_{(m)}}z^{m},\qquad a,b,c,z\in
\mathbb{C}^{d},
\]
where ${v}_{(r)}:=\prod_{i=1}^{d}{(v_i)}_{(r_i)}$ for every
$v,r\in\mathbb{R}^{d}.$

\begin{proposition}
\label{iv}
For every $n\in\mathbb{N}^{d}$ denote $n':=(n_1,\ldots,n_{d-1}).$ A
representation for the connection coefficients in (\ref{lc1}) is
\begin{eqnarray}\label{lau}
c^{*}_m(n)&=&\delta_{mn}\frac{{(|\alpha
|)}_{(|n|)}}{|n|!}\mathit{DM}_{\alpha}(m)\nonumber
\\[-8pt]\\[-8pt]
&&{}\times\sum_{j=0}^{|n|}d_{j}\int_{\Delta
_{(d-1)}}R^{\alpha}_{n'}(t)F_{A}(|\alpha|;-m,-j;\alpha,|\alpha
|;t,1-|t|,1)D_{\alpha}(\mathrm{d}t),\nonumber
\end{eqnarray}
where
\begin{eqnarray}\label{lau2}
d_{j}&:=&\sum_{i=0}^{|n'|}{(-|n'|)}_{(i)}\frac{{(|\alpha
|)}_{(|n'|)}{(|\alpha|+2|n'|)}_{(n_d)}}{i!n_d!}\nonumber
\\[-8pt]\\[-8pt]
&&{}\hspace*{13pt}\times F_A(|\alpha
|;-i,-n_d,-j;|\alpha|,|\alpha|+2i,|\alpha|;1,1,1).\nonumber
\end{eqnarray}
\end{proposition}

 The proof relies on a beautiful representation due to Erd\'{e}lyi
\cite{E38}: for every \mbox{$|a|,|z|\in\mathbb{R}$},
$\alpha,k\in\mathbb{R}^{d}$ and $n\in\mathbb{N}^{d},$
%
\begin{equation}\label{e}
\prod_{j=1}^{d}L_{n_j}^{\alpha_{j}}(k_j|z|)=\sum_{s=0}^{|n|}
\phi_{s}(|a|;\alpha;n;k)L_{s}^{|a|}(|z|),
\end{equation}
where
\[
\phi_{s}(|a|;\alpha;n;k)=F_A(|a|;-n,-s;\alpha,|a|;k,1)\prod
_{j=1}^{d}\frac{{(\alpha_j)}_{(n_j)}}{n_j!}.
\]
The full proof of Proposition \ref{iv} involves tedious algebra that
we omit here as not relevant for the
general purposes of the paper.

\begin{remark}
A simplified representation of $c^{*}_m(n)$ in terms of Hahn
polynomials will be given in Section \ref{lau2sec}.
\end{remark}

\begin{remark}
Note that when $|n'|=0$, $c^{*}_{m}(0,\ldots,0,n_d)=1,$ which agrees
with the known identity
%
\begin{equation}\label{ask}
L_{n}^{\alpha+\beta}(x+y)=\sum_{j=0}^{n}L^{\alpha
}_{j}(x)L^{\beta}_{n-j}(y),
\qquad x,y\in\mathbb{R}
\end{equation}
(see
\cite{AAR99}, formula (6.2.35), page 191), an identity with an obvious
extension to the~$d$-di\-mensional case.
\end{remark}

\begin{remark} \label{scale}It is immediate to verify that the
coefficients $c^{*}_{m}(n)$ also satisfy
%
\begin{equation}L^{|\alpha| }_{|n-n'|}(|\beta^{-1} y|)|\beta^{-1}
y|^{|n'|}R_{n'}^{\alpha}\biggl(\frac{y}{|y|}\biggr)=\sum
_{|m|=|n|}\varphi_m c^{*}_m(n) L_{m}^{\alpha}(|\beta^{-1}
y|), \qquad \beta\in\mathbb{R}_{+}.
\end{equation}
\end{remark}

\subsection{Size-biased multiple Laguerre}\label{sec:sblag}

Let $Y^d=(Y_1,\ldots,Y_d)$ be a collection of independent gamma random
variables, each with parameters $(\theta/d,1),$ $i=1,\ldots,d.$ Let
$\ddot{Y}^d$ be the same vector with the coordinates rearranged in
size-biased random order. The proof of the following corollaries
is, at this point, obvious from Proposition \ref{pr:mlag}.

\begin{corollary}A system of polynomials, orthogonal with respect to
the law of $\ddot{Y}^d$, is given by
%
\begin{equation}\label{premlaggem}
\ddot{L}^{|\theta|,d}
_{(|m|,n')}(y)=L_{|m|}^{|\theta|+2|n'|}({|y|})
({|y|})^{|n'|}\ddot{R}_{n'}^{|\theta|,d}\biggl(\frac
{y}{|y|}\biggr),
\end{equation}
$|m|\in\mathbb{N},\ n'\in\mathbb{N}^{d}\dvtx  |n'|\in\mathbb{N}$, with
$\{\ddot{R}_{n}\}$ as in (\ref{dsbdir}).
\end{corollary}

It is possible to derive an
infinite-dimensional version of $\{L^{\alpha\star}_n\}$,
orthogonal with respect to the law of the size-biased point
process $\ddot{Y}^{\infty},$ obtained by $Y^{\infty}$ of Section
\ref{sec:gempp}. Remember that
$\ddot{X}^{\infty}:=\ddot{Y}^{\infty}/|\ddot{Y}^{\infty}|$ has
$\operatorname{GEM}(|\theta|)$ distribution and it is independent of
$|\ddot{Y}^{\infty}|\stackrel{\mathcal{D}}{=}|Y^{\infty}|$, which
has a gamma$(|\theta|)$ law.

\begin{corollary}
\label{cor:gemlag} Let $\ddot{\gamma}_{|\theta|}$ be the
probability distribution of the size-biased sequence
$\ddot{Y}^{\infty}$ obtained by rearranging in size-biased random
order the sequence $Y^{\infty}$ of points of a~Pois\-son process
with
generating functional (\ref{pdmgf}).
The polynomials defined by
%
\begin{equation}\label{mlaggem}
\ddot{L}^{|\theta|}
_{(|m|,n')}(y)=L_{|m|}^{|\theta|+2|n'|}({|y|})
({|y|})^{|n'|}\ddot{R}_{n'}^{|\theta|}\biggl(\frac
{y}{|y|}\biggr)
\end{equation}
for $|m|\in\mathbb{N}, n'\in\mathbb{N}^{\infty}\dvtx |n'|\in\mathbb
{N}$, with
$\{\ddot{R}_{n}\}$ as in (\ref{gemp}), are the limit, as
$d\rightarrow\infty,$ of the polynomials $\{\ddot{L}^{|\theta
|,d}_{(|m|,n')}\}$ defined by (\ref{premlaggem}) and
form an orthogonal system
with respect to $\ddot{\gamma}_{|\theta|}$.\looseness=-1
\end{corollary}

\section{Multivariate Hahn polynomials}\label{sec:hahn}

\subsection{Hahn polynomials on \{1$,\ldots,$N\}}

As for the Laguerre polynomials, we introduce the discrete Hahn
polynomials on $\{1,\ldots,N\}$ with parameters shifted by 1 to
make the notation consistent with the standard probabilistic
notation in the corresponding weight measure.
The Hahn polynomials, orthogonal on $\mathit{DM}_{\alpha,\beta}(n;N),$ are
defined as the hypergeometric series:
%
\begin{equation}\label{hahn}
h^{\alpha,\beta}_{n}(r;N)={}_{3}F_2\left(
\begin{array}{cc}-n,n+\theta-1,-r\\
\alpha,-N
\end{array}
\ \vline\ \ 1\right),\qquad n=0,1,\ldots,N.
\end{equation}
The orthogonality constants are given by
\[
\frac{1}{{u}^{\alpha,\beta}_{N,n}}:=\sum_{r=0}^{N}[h^{\alpha
,\beta}_{n}(r;N)]^{2}\mathit{DM}_{\alpha,\beta}(n;N)\nonumber
=\frac{1}{\matrix{N\choose
n}}\frac{{(\theta+N)}_{(n)}}{{(\theta)}_{(n-1)}}\frac
{1}{\theta+2n-1}\frac{{(\beta)}_{(n)}}{{(\alpha
)}_{(n)}}.
\]
A special point value is (\cite{KMG61}, formula (1.15))
%
\begin{equation} \label{un}
h^{\alpha,\beta}_{n}(N;N)=(-1)^{n}\frac{{(\beta)}_{(n)}}{{(\alpha
)}_{(n)}}.
\end{equation}
Thus if we consider the normalization
\[
q^{\alpha,\beta}_{n}(r;N):=\frac{h^{\alpha,\beta
}_{n}(r;N)}{h^{\alpha,\beta}_{n}(N;N)},
\]
then the new constant is, from (\ref{un}),
%
\begin{eqnarray}\label{un*}
\frac{1}{{w}^{\alpha,\beta}_{N,n}}&:=&\mathbb{E}[{q}^{\alpha
,\beta}_{n}(R;N)]^{2}\nonumber
\\
&=&\frac{1}{\matrix{N\choose
n}}\frac{{(\theta+N)}_{(n)}}{{(\theta)}_{(n-1)}}\frac
{1}{\theta+2n-1}\frac{{(\alpha)}_{(n)}}{{(\beta
)}_{(n)}}
\\
&=&\biggl[\frac{{(\theta+N)}_{(n)}}{N_{[n]}}\biggr]\frac
{1}{\zeta_n^{\alpha,\beta}},\nonumber
\end{eqnarray}
where $\zeta_n$ is the Jacobi orthogonality constant, given by (\ref
{hn}).

A symmetry relation is
%
\begin{equation}\label{hsym}
q^{\alpha,\beta}_{n}(r;N)=\frac{q^{\beta,\alpha
}_{n}(N-r;N)}{q^{\beta,\alpha}_{n}(0;N)} .
\end{equation}
A well-known relationship is in the limit:
%
\begin{equation}\label{jlim}
\lim_{N\rightarrow\infty}h_n^{\alpha,\beta}(Nz;N)=\widetilde
{R}_{n}^{\alpha-1,\beta-1}(1-2z), \qquad\alpha,\beta>0
\end{equation}
(see \cite{KMG61}), where
$\widetilde{R}^{a,b}_{n}=\widetilde{R}^{a,b}_{n}/\widetilde
{R}^{a,b}_{n}(1)$ are
standardized Jacobi polynomials orthogonal on $[-1,1]$ as defined
in Section \ref{*}. Because of our definition (\ref{01J}),
combining (\ref{symj}), (\ref{hsym}) and (\ref{jlim5.6}) gives the
equivalent limit: For every $n$,
%
\begin{equation}\label{jlim5.6}
\lim_{N\rightarrow\infty}q_n^{\alpha,\beta}(Nz;N)={R}_{n}^{\alpha
,\beta}(z), \qquad \alpha,\beta>0.
\end{equation}
Note that also
%
\begin{equation}\label{limcc}
\lim_{N\rightarrow\infty}{{w}^{\alpha,\beta
}_{N,n}}=\zeta_n^{\alpha,\beta}.
\end{equation}
An inverse relation holds as well, which allows one to derive Hahn
polynomials as a~mixture of Jacobi
polynomials. Denote by $B_x(r;N)=B_{x,1-x}(r,N-r)$ the binomial
distribution.\vadjust{\goodbreak}

\begin{proposition}
\label{pr:uhahn} The functions
%
\begin{eqnarray}\label{uhahn1}
\tilde{q}^{\alpha,\beta}_{n}(r;N)&:=&\int_0^1R_{n}^{\alpha
,\beta
}(x)\frac{B_x(r;N)}{\mathit{DM}_{\alpha,\beta}(r;N)}D_{\alpha,\beta
}(\mathrm{d}x)
\\\label{uhahn2}
&=&\int_0^1R_{n}^{\alpha,\beta}(x)D_{\alpha+r,\beta
+N-r}(\mathrm{d}x),\qquad n=0,1,\ldots,N,
\end{eqnarray}
form the Hahn system of orthogonal polynomials with
$\mathit{DM}_{\alpha,\beta}$ as the weight function, such that
%
\begin{equation}\label{2h}
\tilde{q}^{\alpha,\beta}_{n}(r;N)=\frac{N_{[n]}}{{(\theta
+N)}_{(n)}}
q^{\alpha,\beta}_{n}(r;N).
\end{equation}
\end{proposition}

 The representation (\ref{uhahn2}), in particular, shows
a Bayesian interpretation of Hahn polynomials, as a
\emph{posterior} mixture of Jacobi polynomials evaluated on a
random Bernoulli probability of success $X$, conditionally on
having previously observed $r$ successes out of~$N$ independent
$\operatorname{Bernoulli}(X)$ trials, where $X$ has a $\operatorname{Beta}(\alpha,\beta)$
distribution on $\{0,\ldots,N\}$.

\begin{pf*}{Proof of Proposition \ref{pr:uhahn}}
The integral defined by (\ref{uhahn1}) is a polynomial: Consider
\begin{eqnarray*}
\int_0^1x^{n}(1-x)^{m}\frac{B_x(r;N)}{\mathit{DM}_{\alpha,\beta
}(r;N)}D_{\alpha,\beta}(\mathrm{d}x)&=&\frac{{(\alpha)}_{(n+r)}{(\beta
)}_{(N+m-r)}{(\theta)}_{(N)}}{{(\alpha)}_{(r)}{(\beta
)}_{(N-r)}{(\theta)}_{(N+n+m)}}
\\
&=&\frac{{(\alpha+r)}_{(n)}{(\beta+N-r)}_{(m)}}{{(\theta+N)}_{(n+m)}}.
\end{eqnarray*}
The numerator is a polynomial in $r$ of order $n+m.$ Write
\begin{eqnarray*}
R_{n}^{\alpha,\beta}(x)= \sum_{j=1}^{n}c_jx^{j},
\end{eqnarray*}
then
\begin{eqnarray}\label{L}
\int_0^1
R_{n}^{\alpha,\beta}(x)\frac{B_x(r;N)}{\mathit{DM}_{\alpha,\beta
}(r;N)}D_{\alpha,\beta}(\mathrm{d}x)&=&\sum_{j=1}^{n}\frac{c_j}{{(\theta
+N)}_{(j)}}{(\alpha+r)}_{(j)}\nonumber
\\[-10pt]\\[-10pt]
&=&\sum_{j=1}^{n}\frac
{c_j}{{(\theta+N)}_{(j)}}r_{[j]}+L,\nonumber
\end{eqnarray}
where $L$ is a polynomial in $r$ of order less than $n$. Then
$q_n^{\alpha,\beta}(r)$ is a polynomial of order~$n$ in
$r.$\vadjust{\goodbreak}

To show orthogonality it is sufficient to show that $h_n$ are
orthogonal with respect to polynomials of the basis formed by the
falling factorials $\{r_{[l]}, l=0,1,\ldots\}$. For $l\leq n,$
%
\begin{eqnarray}\label{orh}
&&\sum_{r=0}^{n}\mathit{DM}_{\alpha,\beta}(r;N)r_{[l]}\tilde{q}^{\alpha
,\beta
}_{n}(r;N)\nonumber
\\
&&\quad=\frac{N!}{(N-l)!}\int_0^1x^lR_{n}^{\alpha,\beta
}(x)\Biggl[\sum_{r=0}^{n}\pmatrix{{N-l}\cr{r-l}}x^{l-r}(1-x)^{N-r}
\Biggr]D_{\alpha,\beta}(\mathrm{d}x)
\\
&&\quad=N_{[l]}\int_0^1 x^lR_{n}^{\alpha,\beta}(x)D_{\alpha,\beta}(\mathrm{d}x).\nonumber
\end{eqnarray}
The last integral is non-zero only if $l=n$, which
proves the orthogonality of $q^{\alpha,\beta}_{n}(r;N).$

Now consider that, in $R^{\alpha,\beta}_{n}(x)$, the leading
coefficient $c_{n}$ satisfies
\begin{eqnarray*}
\int_0^1 c_n
x^{n}R^{\alpha,\beta}_{n}(x)D_{\alpha,\beta}(\mathrm{d}x)=\int_0^1
[R^{\alpha,\beta}_{n}(x)]^{2}D_{\alpha,\beta
}(\mathrm{d}x)=\frac{1}{\zeta_n^{\alpha,\beta}};
\end{eqnarray*}\vspace*{-16pt}
\begin{eqnarray}\label{on}
\frac{1}{\omega_{N,n}^{\alpha,\beta}}&=&\sum_{r=0}^{n}\mathit{DM}_{\alpha
,\beta}(r;N)\tilde{q}^{\alpha,\beta}_{n}(r;N)\tilde
{q}^{\alpha,\beta
}_{n}(r;N)\nonumber
\\
&=&
\sum_{r=0}^{n}\mathit{DM}_{\alpha,\beta}(r;N)\Biggl(\sum_{j=0}^{n}\frac
{c_j}{{(\theta+N)}_{(j)}}r_{[j]}\Biggr)\tilde{q}^{\alpha,\beta
}_{n}(r;N)+L'\nonumber
\\[-10pt]\\[-10pt]
&=&N_{[n]}\frac{c_n}{{(\theta+N)}_{(n)}}\int x^{n}R^{\alpha
,\beta}_{n}(x)D_{\alpha,\beta}(\mathrm{d}x)\nonumber
\\
&=&\frac{N_{[n]}}{{(\theta+N)}_{(n)}}\frac{1}{\zeta_n^{\alpha
,\beta}}.\nonumber
\end{eqnarray}
That is,
%
\begin{equation}\label{omwn}
\omega_{N,n}^{\alpha,\beta}=\biggl[\frac{{(\theta
+N)}_{(n)}}{N_{[n]}}\biggr]^{2}w_{N,n}^{\alpha,\beta}
\end{equation}
with $w_{N,n}^{\alpha,\beta}$ as in (\ref{un}), and therefore the
identity (\ref{2h}) follows, completing the proof.
\end{pf*}

\subsection{Multivariate polynomials on the Dirichlet-multinomial distribution}

Multivariate polynomials orthogonal with respect to $\mathit{DM}_{\alpha}$
on the discrete $d$-dimensional simplex were first introduced by
Karlin and McGregor \cite{KMG75} as eigenfunctions of the
birth-and-death process with neutral mutation. Here we derive an
alternative derivation as a~posterior mixture of multivariate
Jacobi polynomials, which extends Proposition \ref{pr:uhahn} to a
multivariate setting.

\begin{proposition}
\label{prp:mh}
For every $\alpha\in\mathbb{R}^{d}$, a system of polynomials,
orthogonal with respect to~$\mathit{DM}_\alpha$, is given by
%
\begin{eqnarray}\label{mh1}
\tilde{q}^{\alpha}_{n}(r;|r|)&=&\int_{\Delta
_{(d-1)}}R_{n}^{\alpha
}(x)\frac{B_{x}(r)}{\mathit{DM}_{\alpha}(r)}D_{\alpha}(\mathrm{d}x)
\\\label{mh2}
&=&\int_{\Delta_{(d-1)}}R_{n}^{\alpha}(x)D_{\alpha+r}(\mathrm{d}x),
\qquad|n|\leq|r|
\\\label{mh3}
&=&\biggl(\frac{\prod
_{j=1}^{d-1}{(A_j+R_j+N_{j+1})}_{(n_{j+1})}}{{(|\alpha
|+|r|)}_{(N_1)}}
\biggr)\prod_{j=1}^{d}\tilde{q}^{\alpha
_j,A_j+2N_j}_{n_j}(r_{j};R_{j-1}-N_j),\quad
\end{eqnarray}
with constant of orthogonality given by
%
\begin{equation}\label{mhc}
\frac{1}{\omega_{n}(\alpha;|r|)}:=\mathbb{E}[\tilde
{q}^{\alpha
}_{n}(R;|r|)]^2=\frac{|r|_{[n]}}{{(|\alpha
|+|r|)}_{(n)}}\frac{1}{\zeta^{\alpha}_{n}}.
\end{equation}
\end{proposition}

\begin{pf}
The identity between (\ref{mh1}) and (\ref{mh2}) is obvious from
Section \ref{sec:conjd} and (\ref{mh3}) follows from Proposition
\ref{pr:uhahn} and some simple algebra. For every
$n\in\mathbb{N}^{d}$,
\begin{eqnarray}\label{mL}
\int_{\Delta_{(d-1)}}x^{n}D_{\alpha+r}(\mathrm{d}x)&=&{\mathit{DM}_{\alpha
+r}(n)}
=\prod_{i=1}^{d-1}\frac{{(\alpha
_i+r_i)}_{(n_i)}{(A_i+R_i)}_{(N_i)}}{{(A_{i-1}+R_{i-1})}_{(N_{i-1})}}\nonumber
\\[-8pt]\\[-8pt]
&=&\frac{\prod_{i=1}^{d}{(\alpha_i+r_i)}_{(n_i)}}{{(|\alpha
|+|r|)}_{(|n|)}}
=\frac{1}{{(|\alpha|+|r|)}_{(|n|)}}\prod
_{i=1}^{d}{r_i}_{[n_i]}+L,\nonumber
\end{eqnarray}
where $L$ is a polynomial in $r$ of order less
than $|n|.$ Therefore $\tilde{q}^{\alpha}_{n}(r;|r|)$ are
polynomials of
order $|n|$ in $r$.

To show that they are orthogonal,
denote
\[
p_l(r):=\prod_{i=1}^{d}(r_i)_{[l_i]}
\]
and consider that, for every $l\in\mathbb{N}^{d}\dvtx |l|\leq|n|,$
%
\begin{eqnarray}\label{morh}
&&\sum_{|m|=|r|}\mathit{DM}_{\alpha}(m;|r|)p_{l}(m)\tilde{q}^{\alpha
}_{n}(m;|r|)\nonumber
\\
&&\quad =\frac{|r|!}{(|r|-|l|)!}\int
x^{l}R_{n}^{\alpha}(x)\biggl(\sum_{|m|=|r|}\pmatrix{{|r-l|}\cr
{m-l}}x^{m-l}\biggr)D_{\alpha}(\mathrm{d}x)
\\
&&\quad=|r|_{[|l|]}\int x^{l}R_{n}^{\alpha}(x)D_{\alpha}(\mathrm{d}x),\nonumber
\end{eqnarray}
which, by orthogonality of $R_{n},$ is non-zero only if $|l|=|n|. $
Since it is always possible to write, for appropriate coefficients
$c_{nm}$
\[
R_{n}^{\alpha}(x)=\sum_{|m|=|n|}c_{nm}x^m+C,
\]
where $C$ is a polynomial of order less than $|n|$ in $x$; then
\[
\tilde{q}^{\alpha}_{s}(r;|r|)=\sum_{|m|=|s|}\frac
{c_{sm}}{(|\alpha
|+|r|)_{(|s|)}}p_m(r)+C'
\]
and by (\ref{morh})
\begin{eqnarray*}
\mathbb{E}[\tilde{q}^{\alpha}_{s}(R;|r|)\tilde
{q}^{\alpha
}_{n}(R;|r|)]&=&\sum_{|k|=|s|}\frac{c_{sk}}{(|\alpha
|+|r|)_{(|s|)}}\mathbb{E}[p_k(R)\tilde{q}^{\alpha
}_{n}(R;|r|)
]+C{''}
\\
&=&|r|_{[|n|]}\sum_{|k|=|r|}\frac{c_{sk}}{(|\alpha
|+|r|)_{(|s|)}}\int x^k R_{n}^{\alpha}(x)D_{\alpha}(\mathrm{d}x)
\\
&=&\frac{|r|_{[|n|]}}{(|\alpha|+|r|)_{(|n|)}}\frac{1}{\zeta
_{n}^{\alpha}}\delta_{sn},\qquad |n|=|r|.
\end{eqnarray*}
\upqed\end{pf}

\begin{remark} Note that the representation (\ref{mh3}) holds also
for negative parameters, so that, if we replace $\alpha$ with $
-\varepsilon$ $(\varepsilon\in\mathbb{R}^d)$ then (\ref{mh3}) is a
representation for polynomials with respect to the hypergeometric
distribution (Section \ref{sec:hg}).
\end{remark}

\subsubsection{Bernstein--B\'{e}zier coefficients of Jacobi polynomials}\label{sec:bb1}

As anticipated in the introduction,
Proposition \ref{prp:mh} gives a probabilistic proof of a recent
result of \cite{W06}, namely that Hahn polynomials are the
Berstein--B\'{e}zier coefficients of the multivariate Jacobi
polynomials. Remember that the Bernstein polynomials, when taken
on the simplex, are essentially multinomial distributions
$B_{x}(n)=\matrix{|n|\choose n}x^{n}$, seen as functions of $x$.

\begin{corollary}\label{cor:bb1}
For every $d\in\mathbb{N},\alpha\in\mathbb{R}^{d},r\in\mathbb{N}^{d}$,
%
\begin{equation}
R_{r}^{\alpha}(x)=\frac{(|\alpha|+|r|)_{(|n|)}}{|r|_{[|n|]}}\sum
_{|m|=|r|}\tilde{q}^{\alpha}_{r}(m;|r|)B_{x}(m),
\label{bb1}
\end{equation}
where $\omega_{r}(|\alpha|;|r|)$ is given by (\ref{mhc}).
\end{corollary}

\begin{pf}
From Proposition \ref{prp:mh},
\[
\mathit{DM}_{\alpha}(m;|r|)\tilde{q}^{\alpha}_{r}(m;|r|)=\mathbb{E}
[B_{X}(m)R_r^{\alpha}(X)]\vspace*{-3pt}
\]
so\vspace*{-6pt}
\[
B_{x}(m)=\mathit{DM}_{\alpha}(m;|m|)\sum_{|n|=0}^{|m|}\zeta_{n}^{\alpha
}\tilde{q}^{\alpha}_{n}(m;|m|)R_{n}(x).\vspace*{-3pt}
\]
Hence\vspace*{-4pt}
%
\begin{eqnarray}
&&\sum_{m}\tilde{q}^{\alpha}_{r}(m;|r|)B_x(m)\nonumber
\\[-3pt]
&&\quad=\sum
_{|n|=0}^{|r|}\zeta
_{n}^{\alpha}\biggl[\sum_{|m|=|r|}\mathit{DM}_{\alpha}(m;|r|)\tilde
{q}^{\alpha
}_{r}(m;|r|)\tilde{q}^{\alpha}_{n}(m;|r|)\biggr]R_{n}^{\alpha
}(x)
\\[-3pt]
&&\quad=\sum_{|n|=0}^{|r|}\frac{\zeta_{n}^{\alpha}}{\omega_{r}(|\alpha
|;|r|)}\delta_{rn}R_{n}^{\alpha}(x)=\frac{|r|_{[|n|]}}{(|\alpha
|+|r|)_{(|n|)}}R_{r}^{\alpha}(x),\nonumber\vspace*{-3pt}
\end{eqnarray}
which completes the proof.\vspace*{-2pt}
\end{pf}

\begin{remark}
By a similar argument it is easy to come back from (\ref{bb1}) to
(\ref{mh1}).\vspace*{-2pt}
\end{remark}

\subsubsection{\texorpdfstring{The connection coefficients of Proposition \ref{iv}}{The connection coefficients of Proposition 4.3}} \label{lau2sec}\vspace*{-1pt}

 Consider again the connection coefficients $c^*_n(m)$ of
Proposition \ref{iv} and their representations~(\ref{lau}) and (\ref
{lau2}). An
alternative representation can be given in terms of multivariate Hahn
polynomials.

\begin{corollary}
Let $c^*_n(m)$ be the connection coefficients between
$L_{n}^{\alpha*}$ and $L_{m}^{\alpha}$, as in Section
\ref{sec:lag}. Then\vspace*{-3pt}
%
\begin{equation}\label{lauhahn}
c^*_n(m)=\delta_{mn}\ b^{|\alpha|}_{|n|,n_d}\mathit{DM}_{\alpha}(m)
\sum_{|r|=0}^{|n|}\frac{{(-m)}_{(r)}}{\prod_{l=1}^{d}r_l!}\tilde
{q}^{\alpha}_{n'}(r;|r|),\vspace*{-4pt}
\end{equation}
where $n'=(n_1,\ldots,n_d-1),$\vspace*{-4pt}
\[
b^{|\alpha|}_{|n|,n_d}=\frac{{(|\alpha|)}_{(|n|)}}{|n|!}
\Biggl[\sum_{j=0}^{|n|}\frac{d_j}{j!|\alpha|_{(j)}}\Biggr]\vspace*{-4pt}
\]
and $d_j$ is as in (\ref{lau2}).
\end{corollary}

\begin{pf}
It is sufficient to use the explicit expression of the Lauricella
function $F_A$ in (\ref{lau}) to see that\vspace*{-3pt}
\begin{eqnarray}
c^{*}_m(n)&=&\delta_{mn}\frac{{(|\alpha
|)}_{(|n|)}}{|n|!}\mathit{DM}_{\alpha}(m)
\Biggl[\sum_{j=0}^{|n|}\frac{d_j}{j!|\alpha|_{(j)}}\Biggr]\sum
_{|r|=0}^{|n|}\frac{{(-m)}_{(r)}}{\prod_{l=1}^{d}r_l!}
\int\frac{\matrix{{|r|}\choose r} t^{r}R_{n'}^{\alpha}(t)}{\mathit{DM}_{\alpha}(r)}D_{\alpha}(\mathrm{d}t)\nonumber
\\[-9pt]\\[-9pt]
&=&\delta_{mn} b^{|\alpha|}_{|n|,n_d}\mathit{DM}_{\alpha}(m)
\sum_{|r|=0}^{|n|}\frac{{(-m)}_{(r)}}{\prod_{l=1}^{d}r_l!}\tilde
{q}^{\alpha}_{n'}(r;|r|).\vspace*{-3pt}\nonumber
\end{eqnarray}
\upqed\end{pf}\nonumber

\subsubsection{Application: The $d$-types linear growth model}\label{sec:kmg}

The multivariate Hahn polynomials were first studied by Karlin and
McGregor \cite{KMG75} to derive the transition density of the
so-called $d$-type neutral Moran model of population genetics. This is,
for any fixed $|r|\in\mathbb{N},$ a stochastic process $(N(t)\dvtx t\geq0)$
living in the discrete simplex $\mathbb{N}_{d,|r|}=\{m\in\mathbb
{N}^{d}\dvtx |m|=|r|\},$ with Dirichlet-multinomial stationary distribution,
and whose generator has Hahn polynomials as eigenfunctions.

Karlin and McGregor's description of such eigenfunctions is
structurally similar to our~(\ref{mh3}), up to some re-scaling and
reordering of the variables.

In the same paper (\cite{KMG75}, formula (6.2)), the functions (rewritten in
our notation)
\[
\psi(m):=\pmatrix{|r|\cr |m|}L_{|r|-|m|}^{|\alpha|+2|m|}(|y|)\tilde
{q}^{\alpha}_{n}(m;|m|),\qquad m\in\mathbb{N}^{d}\dvtx |m|\leq|r|,
|r|\in\mathbb{N},
\]
were introduced to connect the $d$-type Moran model of reproduction to
a $d$-type linear growth model with immigration rates proportional to
$\alpha_1,\ldots,\alpha_d$. The generator of the latter process has
eigenfunctions that are the solution of the recursion
\[
-|y|\psi(m)=\sum_{i=1}^{d}m_i[\psi(m-e_i)-\psi(m)
]+\sum_{i=1}^{d}(m_i+\alpha_i)[\psi(m-e_i)-\psi(m)].
\]
Note that, for every $z\in\mathbb{R}^d$ such that $|z|=|y|$, $\psi
(m)=L_{|r|-|m|,m}^{\alpha}(y)$ is also a solution, hence so is $\psi
(m)=L_{|r|-|m|,m}^{\alpha*}(z).$

 Reconsider now the system $L_{n}^{\alpha*}$ of multiple
Laguerre polynomials. In view of our representation (\ref{mh2}) of
Hahn polynomials, it is easy to write
\[
\psi(m)=\pmatrix{|r|\cr |m|}\frac{\Gamma(|\alpha|)}{\Gamma(\alpha
)}\int_{\mathbb{R}^{d-1}}L_{|r|-|m|,m}^{\alpha*}(y)\frac
{1}{|y|^{d-1}}y^{\alpha-\underline{1}}\,\mathrm{d}y_1\cdots \mathrm{d}y_{d-1},
\]
which is identical to
\[
\psi(m)=\pmatrix{|r|\cr |m|}L_{|r|-|m|}^{|\alpha|+2|m|}(|y|)\int_{\Delta
_{d-1}}R^{\alpha}_{m}(x)D_{\alpha+m}(\mathrm{d}x).
\]
 Our representation in a sense completes Karlin and McGregor's
analysis, in terms of eigenfunctions, of the relationship existing
between the $r$-type linear growth model (product of independent
Laguerre polynomials), the Moran model (multivariate Hahn) and its
scaling limit, the $d$-type Wright--Fisher diffusion (multivariate
Jacobi). In \cite{KMG75} the role of the latter was not very visible.
The representation (\ref{mh2}) shows how to map directly polynomial
eigenfunctions of the scaling limit process (Jacobi ) to polynomial
eigenfunctions of its finite-size dual model (Hahn). In Karlin and
McGregor's work this idea was present only implicitly (see their
formula (3.8) and observation (3.10)), via their use of Laguerre
products. Considering the system $\{L_{|r|-|m|,m}^{\alpha*}\}$ makes
the connection between all the three processes more transparent.


\section{Multivariate Hahn and multiple Meixner polynomials}\label{sec:meix}

The Meixner polynomials on $\{0,1,2,\ldots\},$
defined by
%
\begin{equation}\label{meix}
M_{n}(k;\alpha,p)={}_{2}F_1\left(
\begin{array}{cc}-n,\ -k\\
\alpha
\end{array}
\ \vline\ \
\frac{p-1}{p}\right),\qquad \alpha>0,p\in(0,1),
\end{equation}
are orthogonal with respect to the negative binomial distribution
$ \mathit{NB}_{\alpha,p}.$ The following representation of the Meixner
polynomials comes from the interpretation of $\mathit{NB}_{\alpha,p}$ as
a~gamma mixture of Poisson likelihood (formula (\ref{pgd})).

\begin{proposition}
\label{prp:meix}
For $\alpha\in\mathbb{R}_{+}$ and $p\in(0,1)$, a system of orthogonal
polynomials with the negative binomial $(\alpha,p)$ distribution
as weight measure is given by
%
\begin{eqnarray}\label{meixmix1}
\widetilde{M}^{\alpha,p}_{n}(k)&=&\int_0^{\infty}\frac{\mathit{Po}_{\lambda
}(k)}{\mathit{NB}_{\alpha,p}(k)}L_{n}^{\alpha}\biggl(\lambda
\frac{1-p}{p}\biggr)\gamma_{\alpha,{p/(1-p)}}(\mathrm{d}\lambda)
\\\label{meixmix2}
&=&\int_0^{\infty}L_{n}^{\alpha}\biggl(\lambda\frac{1-p}{p}
\biggr)\gamma_{\alpha+k,p}(\mathrm{d}\lambda),\qquad n=0,1,\ldots,
\end{eqnarray}
where $L_{n}^{\alpha}$ are Laguerre polynomials with parameter $\alpha.$
\end{proposition}

\begin{pf}
For every $n$, consider that
\begin{eqnarray*}
\int_0^{\infty}
\lambda^{n}\gamma_{\alpha+k,p}(\mathrm{d}\lambda)=\int_0^{\infty}\frac
{\lambda^{\alpha+k+n-1}\mathrm{e}^{-{\lambda/p}}}{\Gamma(\alpha
+k)p^{\alpha+k}}\,\mathrm{d}\lambda
={(\alpha+k)}_{(n)}p^n.
\end{eqnarray*}
So every polynomial in $\Lambda$ of order $n$ is mapped to a
polynomial in $k$ of the same order.

To show orthogonality it is, again,
sufficient to consider polynomials in the basis $\{r_{[k]}\dvtx k=0,1,\ldots
\}.$ Let $m\leq n .$
\begin{eqnarray}\label{fmom}
&&\sum_{k=0}^{\infty}\mathit{NB}_{\alpha,p}(k)k_{[m]}\widetilde{M}_{n}^{\alpha
,p}(k)\nonumber
\\
&&\quad =\int_0^{\infty}
L_{n}^{\alpha}\biggl(\lambda\frac{1-p}{p}\biggr)\Biggl\{\sum
_{k=0}^{\infty}\frac{{(\alpha)}_{(k)}}{k!}p^{k}(1-p)^{\alpha
}k_{[m]}\frac{\lambda^{\alpha+k-1}\mathrm{e}^{-{\lambda/p}}}{\Gamma
(\alpha+k)p^{\alpha+k}}\Biggr\}\,\mathrm{d}\lambda\nonumber
\\[-8pt]\\[-8pt]
&&\quad =\int_0^{\infty} L_{n}^{\alpha}\biggl(\lambda\frac{1-p}{p}
\biggr)\Biggl\{\sum_{k=0}^{\infty}k_{[m]}\mathit{Po}_{\lambda}(k)\Biggr\}\gamma
_{\alpha,{p/(1-p)}}(\mathrm{d}\lambda)\nonumber
\\
&&\quad =\int_0^{\infty}
L_{n}^{\alpha}\biggl(\lambda\frac{1-p}{p}\biggr)\lambda^{m}\gamma
_{\alpha,{p/(1-p)}}(\mathrm{d}\lambda),\nonumber
\end{eqnarray}
where the last line comes from the fact that, if $K$ is a
Poisson$(\lambda)$ random variable, then
\[
\mathbb{E}_{\lambda}\bigl(K_{[n]}\bigr)=\lambda^{n},\qquad n=0,1,2,\ldots.
\]
Now, consider the change of measure induced by
\[
z:=\lambda\frac{1-p}{p}.
\]
The last line of (\ref{fmom}) reads
\[
\biggl(\frac{p}{1-p}\biggr)^{m}\int_0^{\infty} L_{n}^{\alpha
}(z)z^m\gamma_{\alpha,1}(\mathrm{d}z).
\]
The integral vanishes for every $m<n,$ and therefore the
orthogonality is proved.
\end{pf}

 From property (2) of the negative binomial distribution
(Section \ref{snbsum}), by using Propositions~\ref{prp:meix}, \ref
{prp:mh} and \ref{iv}, and Remark \ref{scale}, it is possible to
find the following alternative systems of multivariate Meixner
polynomials, orthogonal with respect to $\mathit{NB}_{\alpha,p}^{d}(r).$

\begin{proposition} \label{prp:mmeix} Let $\alpha\in\mathbb
{R}_{+}^{d}$ and $p\in(0,1).$
\begin{longlist}[(ii)]
\item[(i)] Two systems of multivariate orthogonal polynomials with weight
measure $\mathit{NB}_{\alpha,p}^{d}(r)$ are:
%
\begin{equation}\label{malpha}
\widetilde{M}_{n}^{\alpha,p}(r)=\prod_{i=1}^{d}\tilde
{M}_{n_i}^{\alpha
_i,p}(r_i),\qquad n\in\mathbb{N}^{d},
\end{equation}
and
%
\begin{equation}\label{mm}
{}^{*}\widetilde{M}_{n}^{\alpha,p}(r)=(1-p)^{|n'|}\tilde
{M}^{|\alpha
|+2|n'|,p}_{n_d}(|r|-|n'|) {(|\alpha+r|)}_{(|n'|)}\ \tilde
{q}_{n'}^{\alpha}(r;|r|),\qquad n\in\mathbb{N}^{d},
\end{equation}
where $n'=(n_1,\ldots,n_d-1)$, $\{M_{n_i}^{\alpha_i,p}\}$ are Meixner
polynomials as in Proposition \ref{prp:meix} and~$\tilde{q}_{\alpha}$ are
multivariate Hahn polynomials defined by Proposition \ref{prp:mh}.

\item[(ii)] A representation for these polynomials is:
%
\begin{eqnarray}\label{mm11}
\widetilde{M}_{n}^{\alpha,p}(r)&=&\int_{\mathbb{R}_{+}^{d}}\frac
{\mathit{Po}_{\lambda
}^{d}(r)}{\mathit{NB}_{\alpha,p}^{d}(r)}L_{n}^{\alpha}\biggl(\lambda\frac
{1-p}{p}\biggr)\gamma_{\alpha,{p/(1-p)}}^{d}(\mathrm{d}\lambda)
\\\label{mm12}
&=&\int_{\mathbb{R}_{+}^{d}}L_{n}^{\alpha}\biggl(\lambda\frac
{1-p}{p}\biggr)\gamma_{\alpha+r,p}^{d}(\mathrm{d}\lambda)
\end{eqnarray}
and
%
\begin{eqnarray}\label{mm21}
{}^{*}\widetilde{M}_{n}^{\alpha,p}(r)&=&\int_{\mathbb
{R}_{+}^{d}}\frac
{\mathit{Po}_{\lambda}^{d}(r)}{\mathit{NB}_{\alpha,p}^{d}(r)}L_{n}^{\alpha*}
\biggl(\lambda\frac{1-p}{p}\biggr)\gamma_{\alpha,
{p/(1-p)}}^{d}(\mathrm{d}\lambda)
\\\label{mm22}
&=&\int_{\mathbb{R}_{+}^{d}}L_{n}^{\alpha*}\biggl(\lambda\frac
{1-p}{p}\biggr)\gamma_{\alpha+r,p}^{d}(\mathrm{d}\lambda),
\end{eqnarray}
where $\{L_{n}^{\alpha}\}$ and $\{L_{n}^{\alpha*}\}$ are given by
(\ref{mlag1}) and (\ref{mlag}), and
\[
\gamma_{\alpha,\beta}^{d}(\mathrm{d}z):=\prod_{i=1}^{d}\gamma_{\alpha
_i,\beta}(\mathrm{d}z_i),\qquad \beta\in\mathbb{R},z\in\mathbb{R}^{d}.
\]

\item[(iii)] The connection coefficients between $\{\widetilde{M}_{n}^{\alpha
,p}\}$
and ${}^{*}\widetilde{M}_{n}^{\alpha,p}$ are given by
%
\begin{equation}\label{ccmeix}
\mathbb{E}[{}^{*}\widetilde{M}_{n}^{\alpha,p}(R)\tilde
{M}_{m}^{\alpha
,p}(R)]={c}^{*}_{m}(n),
\end{equation}
where ${c}^{*}_{m}(n)$ are as in (\ref{lau}) or (\ref{lauhahn}).
\end{longlist}
\end{proposition}

\begin{pf}
(\ref{malpha}) is trivial and (\ref{mm11}) and (\ref{mm12}) follow from
(\ref{meixmix1}) and (\ref{meixmix2}).

Now let us first prove (\ref{mm21}) and (\ref{mm22}). For every $z\in
\mathbb{R}_{+}^{d},$ denote $x=z/|z|.$ Consider that
\begin{eqnarray*}
\gamma_{\alpha, \beta}(\mathrm{d}z)=\gamma_{|\alpha|,\beta}(\mathrm{d}|z|)D_{\alpha}(\mathrm{d}x)
\end{eqnarray*}
and that
\[
\mathit{Po}_{z}^{d}(r)=\mathit{Po}_{|z|}(|r|)L_{x}(r).
\]
Combining this with (\ref{nbsum}),
%
\begin{eqnarray}\label{decmm}
&&\int_{\mathbb{R}_{+}^{d}}\frac{\mathit{Po}_{\lambda}^{d}(r)}
{\mathit{NB}_{\alpha,p}^{d}(r)}L_{n}^{\alpha*}\biggl(\lambda\frac
{1-p}{p}\biggr)\gamma_{\alpha,{p/(1-p)}}^{d}(\mathrm{d}\lambda)\nonumber
\\
&&\quad =\biggl(\int_{\mathbb{R}_{+}}\frac{\mathit{Po}_{|\lambda
|}(|r|)}{\mathit{NB}_{|\alpha
|,p}(|r|)}L_{n_d}^{|\alpha|+2|n'|}\biggl(|\lambda|\frac
{1-p}{p}\biggr)\biggl[|\lambda|\frac{1-p}{p}\biggr]^{|n'|}\gamma
_{|\alpha|,{p/(1-p)}}(\mathrm{d}|\lambda|)\biggr)\quad
\\
&&\qquad{}\times\biggl(\int_{\Delta_{(d-1)}}\frac
{L_{x}(r)}{\mathit{DM}_{\alpha}(r,|r|)}R_{n'}^{\alpha}(x)D_{\alpha}(\mathrm{d}x)\biggr).
\nonumber
\end{eqnarray}
From Proposition \ref{prp:mh}, the last integral in (\ref{decmm})
is equal to $\tilde{q}_{n'}^{\alpha}(r;|r|).$

The first integral
can be rewritten as
%
\begin{eqnarray}\label{bitmm}
&&\int_{\mathbb{R}_{+}}L_{n_d}^{|\alpha|+2|n'|}\biggl(|\lambda
|\frac
{1-p}{p}\biggr)
\biggl[|\lambda|\frac{1-p}{p}\biggr]^{|n'|}\gamma_{|\alpha
|+|r|,{p/(1-p)}}(\mathrm{d}|\lambda|)\nonumber
\\
&&\quad = (1-p)^{|n'|}{(|\alpha+r|)}_{(|n'|)}\int_{\mathbb{R}_{+}}
L_{n_d}^{|\alpha|+2|n'|}\biggl(|\lambda|\frac{1-p}{p}\biggr)\frac
{|\lambda|^{|\alpha+r+n'|}\mathrm{e}^{-{|\lambda|/p}}}{\Gamma(|\alpha
+r+n'|)p^{|\alpha+r+n'|}}\,\mathrm{d}|\lambda|\qquad\quad
\\
&&\quad=(1-p)^{|n'|}{(|\alpha+r|)}_{(|n'|)}\widetilde{M}_{n_d}^{\alpha
+2|n'|}(|r|-|n'|).\nonumber
\end{eqnarray}
The last line in (\ref{bitmm}) is obtained from (\ref{meixmix2}) by
rewriting $|n'|=2|n'|-|n'|$ in the mixing measure. Thus the identities
(\ref{mm21}) and (\ref{mm22}) are proved.

To prove part (iii), simply use (\ref{lc1}) with coefficients
given by Proposition \ref{iv} to see that
(\ref{mm11}) and (\ref{mm12}) and (\ref{mm21}) and (\ref{mm22}) imply
\begin{eqnarray*}
{}^{*}\widetilde{M}_{n}^{\alpha,p}(r)&=&\mathbb{E}_{\alpha
+r,p}
\biggl[L_{n}^{\alpha*}\biggl(\lambda\frac{1-p}{p}\biggr)\biggr]
=\mathbb{E}_{\alpha+r,p}\biggl[\sum
_{|m|=|n|}c^{*}_m(n)L_{m}^{\alpha
}\biggl(\lambda\frac{1-p}{p}\biggr)\biggr]
\\
&=&\sum_{|m|=|n|}c^{*}_m(n)\mathbb{E}_{\alpha+r,p}
\biggl[L_{m}^{\alpha
}\biggl(\lambda\frac{1-p}{p}\biggr)\biggr]
=\sum_{|m|=|n|}c^{*}_m(n)\widetilde{M}_{m}^{\alpha,p}(r).
\end{eqnarray*}
This is equivalent to (\ref{ccmeix}) because of the orthogonality
of $\widetilde{M}_{m}^{\alpha,p}(R).$

But (\ref{ccmeix}) also implies that $\{
{}^{*}\widetilde{M}_{n}^{\alpha,p}(r)\}$ is an orthogonal system with
$\mathit{NB}_{\alpha,p}^{d}$ as weight measure since, for every polynomial
$r_{[l]}$ of degree $|l|\leq|n|,$
\begin{eqnarray*}
\sum_{r\in\mathbb{N}^{d}}\mathit{NB}_{\alpha,p}^{d}(r){}^{*}\tilde
{M}_{n}^{\alpha
,p}(r)r_{[l]}=\sum_{|m|=|n|}c^{*}_m(n)\biggl(\sum_{r\in\mathbb
{N}^{d}}\mathit{NB}_{\alpha,p}^{d}(r)\widetilde{M}_{m}^{\alpha
,p}(r)r_{[l]}
\biggr).
\end{eqnarray*}
The term between brackets is non-zero only for $|l|=|m|=|n|,$
which implies orthogonality, so the proof of the proposition is
now complete.
\end{pf}

\subsection{The Bernstein--B\'{e}zier coefficients of the multiple Laguerre polynomials}\label{sec:bb2}

The representation of Meixner polynomials given in Proposition
\ref{prp:mmeix} leads us, not surprisingly, to interpret these as the
Bernstein--B\'{e}zier coefficients of the multiple Laguerre
polynomials (for any choice of basis), up to proportionality
constants. Note that, for products of Poisson distributions we can
write
%
\begin{equation}\label{pober}
\mathit{Po}_{\lambda}^{d}(r)=\prod_{i=1}^{d}\frac{\mathrm{e}^{-\lambda_i}\lambda
_i^{r_{i}}}{r_{i}!}=\frac{\mathrm{e}^{-|\lambda|}}{|\lambda|!}B_{\lambda}(r).
\end{equation}
 To simplify the notation, let $(L_{m},M_n)$ denote
either $(L^{\alpha}_{m},\widetilde{M}^{\alpha,p}_m)$ or
$(L^{\alpha*}_{m},{}^*\widetilde{M}^{\alpha,p}_m),$ for some
$\alpha\in\mathbb{R}^{d}$ and $p\in(0,1).$ Let $\varphi_n$ be either
as in (\ref{mlagcp}) or as in (\ref{mlconst}), consistently with
the choice of $L_n$, and set
$\rho_{r}(\alpha,p)^{-1}:=E[M_{r}^2].$\vspace*{-2pt}

\begin{corollary}\vspace*{-3pt}
%
\begin{equation}\label{bb2}
L_{r}\biggl(\lambda\frac{1-p}{p}\biggr)=\frac{\rho_{r}(\alpha
,p)}{\varphi_r}\frac{\mathrm{e}^{-|\lambda|}}{|\lambda|!}\sum
_{m}M_{r}(m)B_{\lambda}(m).
\end{equation}
\end{corollary}

\begin{pf}
The proof is along the same lines as for Corollary \ref{cor:bb1}.
From (\ref{mm11})--(\ref{mm21}),\vspace*{-2pt}
\[
\mathbb{E}\biggl[L_{n}\biggl(Y\frac{1-p}{p}
\biggr)\mathit{Po}_{Y}^{d}(m)
\biggr]=M_{n}(m)\mathit{NB}_{\alpha,p}^{d}(m), \qquad n,m\in\mathbb{N}^{d}.\vspace*{-2pt}
\]
Then from (\ref{pober}),\vspace*{-2pt}
\[
B_{\lambda}(m)=|\lambda|!\mathrm{e}^{|\lambda|}\mathit{NB}_{\alpha,p}^{d}(m)
\sum_{n}\varphi_nM_{n}(m)L_{n}\biggl(Y\frac{1-p}{p}\biggr).\vadjust{\goodbreak}
\]

\noindent
So
for every $r\in\mathbb{N}^{d}$
\begin{eqnarray*}
\sum_{m}M_{r}(m)B_{\lambda}(m)&=&|\lambda|!\mathrm{e}^{|\lambda|}\sum
_{n}\varphi_n\biggl[\sum_{m}\mathit{NB}_{\alpha
,p}^{d}(m)M_{n}(m)M_{r}(m)\biggr]L_{n}\biggl(Y\frac{1-p}{p}
\biggr)\nonumber
\\
&=&|\lambda|!\mathrm{e}^{|\lambda|}\sum_{n}L_{n}\biggl(Y\frac{1-p}{p}
\biggr)\frac{\varphi_n}{\rho_{r}(\alpha,p)}\delta_{nr}\nonumber
\\
&=&\frac{|\lambda|!\mathrm{e}^{|\lambda|}\varphi_r}{\rho_{r}(\alpha
,p)}L_{r}\biggl(Y\frac{1-p}{p}\biggr),
\end{eqnarray*}
and the proof is complete.
\end{pf}

\section*{Acknowledgement}
Dario Span\`{o}'s research is partly supported by CRiSM, an
EPSRC-funded grant.

\printhistory

\end{document}